\documentclass[onefignum,onetabnum]{siamart190516}

\usepackage{graphicx}
\usepackage{tabularx}
\usepackage{multirow}

\usepackage{latexsym}
\usepackage{amsmath,amstext,amssymb,bbm,amsbsy}
\usepackage{mathtools}
\usepackage{bm}
\usepackage{caption}
\usepackage{subcaption}
\usepackage{xcolor}

\usepackage{alphalph}

\usepackage[multiple]{footmisc}

\newtheorem{remark}{Remark}
\newtheorem{example}{Example}

\newcommand{\R}{\mathbb{R}}
\newcommand{\Rn}{\mathbb{R}^n}
\newcommand{\ba}{\boldsymbol{a}}
\newcommand{\bb}{\boldsymbol{b}}
\newcommand{\bx}{\boldsymbol{x}}
\newcommand{\bv}{\boldsymbol{v}}
\newcommand{\bp}{\boldsymbol{p}}
\newcommand{\by}{\boldsymbol{y}}
\newcommand{\bq}{\boldsymbol{q}}
\newcommand{\bu}{\boldsymbol{u}}
\newcommand{\bz}{\boldsymbol{z}}
\newcommand{\bw}{\boldsymbol{w}}
\newcommand{\bzero}{\boldsymbol{0}}
\newcommand{\bmu} {\boldsymbol{\mu}}
\newcommand{\blambda} {\boldsymbol{\lambda}}
\newcommand{\norm}[1]{\left\lVert#1\right\rVert}

\title{SympOCnet: Solving optimal control problems with applications to high-dimensional multi-agent path planning problems\thanks{Submitted to the editors \today.
\funding{This work was supported by
OSD/AFOSR MURI grant FA9550-20-1-0358.}}}
\author{Tingwei Meng\footnotemark[2]\ \footnotemark[3] \and
Zhen Zhang\footnotemark[2]\ \footnotemark[3] \and
J\'er\^ome Darbon\footnotemark[2]\ \footnotemark[4] \and
George Em Karniadakis\footnotemark[2]\ \footnotemark[5]
}

\begin{document}
\maketitle
\renewcommand{\thefootnote}{\fnsymbol{footnote}}
\footnotetext[2]{Division of Applied Mathematics, Brown University, Providence, RI 02912, USA (tingwei\_meng@brown.edu, zhen\_zhang1@brown.edu, Jerome\_Darbon@brown.edu, George\_Karniadakis@Brown.edu).}
\footnotetext[3]{Tingwei Meng and Zhen Zhang contributed equally to this work.}
\footnotetext[4]{Corresponding author.}
\footnotetext[5]{School of Engineering, Brown University, Providence, RI 02912, USA.}
\renewcommand{\thefootnote}{\arabic{footnote}}

\begin{abstract}
    Solving high-dimensional optimal control problems in real-time is an important but challenging problem, with applications to multi-agent path planning problems, which have drawn increased attention given the growing popularity of drones in recent years. In this paper, we propose a novel neural network method called SympOCnet that applies the Symplectic network to solve high-dimensional optimal control problems with state constraints. We present several numerical results on path planning problems in two-dimensional and three-dimensional spaces. Specifically, we demonstrate that our SympOCnet can solve a problem with more than $500$ dimensions in 1.5 hours on a single GPU,
    which shows the effectiveness and efficiency of SympOCnet. The proposed method is scalable and has the potential to solve truly high-dimensional path planning problems in real-time.
\end{abstract}

\begin{keywords}
  deep neural networks, optimal control, path planning, physics-informed learning
\end{keywords}

\begin{AMS}
  49M99, 68T07
\end{AMS}

\section{Introduction}
Optimal control methods are widely applied in many practical problems, which include path planning~\cite{Coupechoux2019Optimal,Rucco2018Optimal,Hofer2016Application,Delahaye2014Mathematical,Parzani2017HJB,Lee2021Hopf,Makkapati2021Desensitized}, humanoid robot control~\cite{Khoury2013Optimal,Feng2014Optimization,kuindersma2016optimization,Fujiwara2007optimal,fallon2015architecture,denk2001synthesis}, and robot manipulator control~\cite{lewis2004robot,Jin2018Robot,Kim2000intelligent,Lin1998optimal,Chen2017Reachability}.
With the increasing impact of drones in everyday tasks as well as specialized military tasks, the path planning problem of multiple drones has become increasingly important. However, it is challenging to model this problem using the optimal control formulation because the dimensionality is usually very high, causing difficulty in applying traditional numerical methods, such as the very accurate pseudospectral (PS) method~\cite{fahroo2001costate,fahroo2002direct,Rao2009Survey,Ross2012Review,garg2011advances}, to solve it. 
For instance, for a multi-agent path planning problem with $M$ drones, the motion of each drone has to be described by a state variable in $\R^m$, therefore the dimension of the corresponding optimal control problem will be $Mm$, which is huge when the number of drones $M$ is large. Computational complexity of traditional numerical methods may scale exponentially with respect to the dimension, which makes the high-dimensional problems intractable. This is an open issue called ``curse-of-dimensionality'' (CoD)~\cite{bellman1961adaptive}. Hence, new approaches are required to tackle high-dimensional optimal control problems efficiently.

In literature, there are some algorithms for solving high-dimensional optimal control problems (or the corresponding Hamilton-Jacobi PDEs), which include optimization methods~\cite{darbon2015convex,Darbon2016Algorithms,darbon2019decomposition,darbon2021hamilton,Chen2021Lax,Chen2021Hopf,yegorov2017perspectives,Lee2021Computationally,Kirchner2020HJ}, 
max-plus methods~\cite{akian2006max,akian2008max,dower2015maxconference,Fleming2000Max,gaubert2011curse,McEneaney2006maxplus,McEneaney2007COD,mceneaney2008curse,mceneaney2009convergence}, 
tensor decomposition techniques \cite{dolgov2019tensor,horowitz2014linear,todorov2009efficient}, sparse grids \cite{bokanowski2013adaptive,garcke2017suboptimal,kang2017mitigating}, 
polynomial approximation \cite{kalise2019robust,kalise2018polynomial}, 
model order reduction \cite{alla2017error,kunisch2004hjb}, optimistic planning~\cite{Bokanowski2021Optimistic},
dynamic programming and reinforcement learning \cite{Chen2015Exact,Chen2015Safe,alla2019efficient,bertsekas2019reinforcement,zhou2021actor}, as well as methods based on neural networks \cite{bachouch2018deep,bansal2020deepreach,Djeridane2006Neural,jiang2016using,Han2018Solving,hure2018deep,hure2019some,lambrianides2019new,Niarchos2006Neural,reisinger2019rectified,royo2016recursive,Sirignano2018DGM,Li2020generating,darbon2020overcoming,Darbon2021Neural,darbon2021neuralcontrol,nakamurazimmerer2021adaptive,NakamuraZimmerer2021QRnet,jin2020learning,jin2020sympnets,onken2021neural}. 
However, path planning problems are in general hard to solve and cannot be solved directly by applying most of the aforementioned algorithms. One difficulty comes from the state constraints given by the complicated environment. 
Therefore, solving high-dimensional optimal control problems with complicated state constraints is an important yet challenging open problem.
The difficulty of high dimensionality can be avoided by solving the sequential path planning problem instead of the simultaneous planning problem~\cite{Chen2015Exact,Chen2015Safe,Robinson2018Efficient}. The main idea of sequential planning is to assign for different drones different levels of priorities and then sequentially solve the low-dimensional path planning problem involving each drone according to its priority.
However, sequential planning only provides a feasible solution to the original simultaneous planning problem, which may not be optimal. Hence, how to solve high-dimensional simultaneous planning problems with complicated state constraints is still an unsolved problem.

Recently, neural networks have achieved great success in solving high-dimensional PDEs and inverse problems and bear some promise in overcoming the CoD. We can make further progress if we take advantage of the knowledge about the underlying data generation process and encode such information in the neural network architectures and the loss functions \cite{jin2020learning,zhang2021gfinns,greydanus2019hamiltonian,meng2020composite}. As a specific example, the physics-informed neural network (PINN), which was first introduced in \cite{raissi2019physics}, encodes physical information directly into the loss function of neural networks using the residuals of ODEs/PDEs and automatic differentiation.
In this work, we leverage PINNs and propose a novel method using a Symplectic optimal control network called SympOCnet to solve high-dimensional optimal control problems. SympOCnet is designed based on the Symplectic network (SympNet), which is a neural network architecture proposed in~\cite{jin2020sympnets} and is able to approximate arbitrary symplectic maps. In the original paper, SympNet is used to solve the inverse parameter identification problem by acting as a surrogate model for the data generation procedure. In this paper, we present a novel way to utilize SympNet for solving Hamiltonian ODEs in a forward manner.
By applying the SympOCnet, we encode our knowledge about optimal control problems and Hamilton-Jacobi theory in the neural network. We then apply SympOCnet on multi-agent path planning problems and show its effectiveness and efficiency in solving high-dimensional problems. Specifically, we can solve a path planning problem involving $256$ drones, which corresponds to a $512$-dimensional optimal control problem.

The organization of this paper is given as follows. In Section~\ref{sec:bkgd}, we briefly summarize some preliminary background which will be used in the paper, including Hamiltonian systems and symplectic maps in Section~\ref{sec:bkgd_HODE_symplectic} and SympNet in Section~\ref{sec:bkgd_sympnet}. In Section~\ref{sec:optctrl_unconstrained}, we present our proposed SympOCnet method for the optimal control problems without state constraints, which serves as a building block for the full version of the SympOCnet method. 
The full version is then introduced in Section~\ref{sec:optctrl_constrained}, which is proposed for solving the optimal control problems with state constraints. To be specific, we propose four loss functions, which are compared in the numerical experiments. Our numerical experiments on multi-agent path planning problems are presented in Section~\ref{sec:numerical_results}. The first example in Section~\ref{sec:test1_comparison} is used to compare the four loss functions. The second example in Section~\ref{sec:test2_uncertain} demonstrates the generalizability of our method when only partial data is available in the training process. Then, in Sections~\ref{sec:test3_multidrones} and~\ref{sec:test4_swarm}, we demonstrate the effectiveness and efficiency of SympOCnet on high-dimensional problems, with dimensions ranging from $64$ to $512$.
A summary is provided in Section~\ref{sec:conclusion}.

\section{Preliminary background}\label{sec:bkgd}
In this paper, we solve optimal control problems with SympOCnet, a neural network with SympNet architecture and loss function given by the corresponding Hamiltonian ODEs.
This section provides some mathematical materials used in the remainder of the paper. 
In Section~\ref{sec:bkgd_HODE_symplectic}, we provide a brief summary about symplectic maps, Hamiltonian ODEs and their relations. Then, in Section~\ref{sec:bkgd_sympnet}, we review the SympNet architectures. For more details regarding the theory and applications of symplectic methods, we refer the readers to \cite{hairer2006geometric}.

\subsection{Hamiltonian systems and symplectic maps} \label{sec:bkgd_HODE_symplectic}

\begin{definition}[Symplectic maps]\label{def:symp_map}
Let $U$ be an open set in $\R^{2n}$. A differentiable map $\phi:U\rightarrow\R^{2n}$ is called symplectic if the Jacobian matrix $\nabla \phi$ satisfies
\begin{equation}\label{eq:symp_map}
    \nabla \phi^T(\bz) J \nabla \phi(\bz) = J \quad\forall \bz\in\R^{2n}, 
\end{equation}
where $J$ is a matrix with $2n$ rows and $2n$ columns defined by
\begin{equation}\label{eqt:bkgd_def_J}
    J:=\begin{pmatrix} \bm{0} & I_{n} \\ -I_{n} & \bm{0} \end{pmatrix},
\end{equation}
and $I_n$ denotes the identity matrix with $n$ rows and $n$ columns.
\end{definition}
The Hamiltonian ODE system is a dynamical system taking the form 
\begin{equation*}
    \dot{\bz}(s) = J\nabla H(\bz(s)),
\end{equation*}
where $J$ is the matrix defined in~\eqref{eqt:bkgd_def_J} and $H:U\rightarrow \R$ is a function called Hamiltonian.
The Hamiltonian systems and symplectic maps are highly related to each other. To be specific, the Hamiltonian structure is preserved under the change of variable using any symplectic map. This result is stated in the following theorem. 
\begin{theorem}\cite[Theorem~2.8 on p.~187]{hairer2006geometric}\label{thm:bkgd_symplectic_HODE}
Let $U$ and $V$ be two open sets in $\R^{2n}$. Let $\phi: U \rightarrow V$ be a change of coordinates such that $\phi$ and $\phi^{-1}$ are continuously differentiable functions. If $\phi$ is symplectic, the Hamiltonian ODE system $\dot{\bz}(s) = J\nabla H(\bz(s))$ can be written in the new variable $\bw = \phi(\bz)$ as
\begin{equation}\label{eq:change_var}
    \dot{\bw}(s) = J\nabla \tilde{H}(\bw(s)),
\end{equation}
where the new Hamiltonian $\tilde{H}$ is defined by 
\begin{equation}\label{eqt:thm1_def_tildeH}
    \tilde{H}(\bw) = H(\bz) = H(\phi^{-1}(\bw)) \quad\forall \bw\in V.
\end{equation}
Conversely, if $\phi$ transforms every Hamiltonian system to another Hamiltonian system by \eqref{eq:change_var} and~\eqref{eqt:thm1_def_tildeH}, then $\phi$ is symplectic.
\end{theorem}
Theorem~\ref{thm:bkgd_symplectic_HODE} indicates that a symplectic map can transform a Hamiltonian ODE system to another system which is again Hamiltonian and potentially in a simpler form. 
This is the starting point of our proposed method. It is well known from literature that the optimal trajectory of an optimal control problem is related to a Hamiltonian ODE system under certain assumptions. Such systems may be solved more easily if written in the right coordinates. We propose to learn the corresponding coordinate transformation through a parameterized family of symplectic maps $\phi_{\theta}$, where $\theta$ represents the unknown parameters to be learned, and solve the Hamiltonian ODEs with new coordinates. The solution to the original problem can be obtained by mapping the trajectory back to original phase space through $\varphi_{\theta} = \phi_{\theta}^{-1}$.

\subsection{Symplectic networks (SympNets)}\label{sec:bkgd_sympnet}
SympNet is a neural network architecture proposed in~\cite{jin2020sympnets} to approximate symplectic transformations. 
There are different kinds of SympNet architectures. In this paper, we use the G-SympNet as the class of our symplectic transformation $\varphi_{\theta}$. For other architectures, we refer the readers to~\cite{jin2020sympnets}.

Define a function $\hat{\sigma}_{K,\ba,\bb}\colon \R^n\to\Rn$ by $\hat{\sigma}_{K,\ba,\bb}(\bx):=K^T(\ba\odot\sigma(K\bx+\bb))$ for any $\bx\in\R^n$, where $\sigma$ is the activation function, and $\odot$ denotes the componentwise multiplication. Here $\ba,\bb$ are vectors in $\R^l$, $K$ is a matrix with $l$ rows and $n$ columns.
Any G-SympNet is an alternating composition of the following two parameterized functions:
 \begin{equation}\label{eq:symp_net}
    \begin{split}
    &\mathcal{G}_{up}\begin{pmatrix} \bx \\ \bp \end{pmatrix}=\begin{pmatrix} \bx\\\bp+\hat{\sigma}_{K,\ba,\bb}(\bx) \end{pmatrix} \quad \forall \bx,\bp\in\R^n,\\ &\mathcal{G}_{low}\begin{pmatrix} \bx \\ \bp \end{pmatrix}=\begin{pmatrix}  \hat{\sigma}_{K,\ba,\bb}(\bp)+\bx \\ \bp \end{pmatrix}\quad \forall \bx,\bp\in\R^n,    
    \end{split}
    \end{equation}
    where the learnable parameters are the matrix $K$ and the vectors $\ba,\bb\in\R^l$.
    The dimension $l$ (which is the dimension of $\ba,\bb$ as well as the number of rows in $K$) is a positive integer that can be tuned, called the width of SympNet. 
In \cite{jin2020sympnets}, it is proven that G-SympNets are universal approximators within the family of symplectic maps.
Note that it is easy to obtain the inverse map of a G-SympNet, since we have explicit formulas for $\mathcal{G}_{up}^{-1}$ and $\mathcal{G}_{low}^{-1}$ given as follows
\begin{equation}\label{eq:symp_inverse}
    \mathcal{G}_{up}^{-1}\begin{pmatrix} \bx\\\bp \end{pmatrix}=\begin{pmatrix}
    \bx \\ \bp-\hat{\sigma}_{K,\ba,\bb}(\bx)  \end{pmatrix},\quad \mathcal{G}_{low}^{-1}\begin{pmatrix}  \bx\\\bp \end{pmatrix}=\begin{pmatrix}  \bx-\hat{\sigma}_{K,\ba,\bb}(\bp) \\ \bp \end{pmatrix}.
    \end{equation}
\section{SympOCnet for optimal control problems without state constraints} \label{sec:optctrl_unconstrained}

In this section, we consider the following optimal control problem without state constraints
\begin{equation} \label{eqt:optctrl_unconstrained}
\begin{split}
    &\min\left\{\int_0^T
    \ell(\bx(s),\bv(s)) ds \colon \bx(0)=\bx_0, \bx(T)=\bx_T, \dot{\bx}(s) = \bv(s)\in U\, \forall s\in(0,T) \right\}\\
    =\ &\min\left\{\int_0^T
    \ell(\bx(s),\dot{\bx}(s)) ds \colon \bx(0)=\bx_0, \bx(T)=\bx_T, \dot{\bx}(s) \in U\, \forall s\in(0,T) \right\},
\end{split}
\end{equation}
where $\bv\colon [0,T]\to U$ is the control function taking values in the control set $U\subseteq \Rn$, $\bx\colon [0,T]\to \Rn$ is the corresponding trajectory, and $\ell\colon \Rn\times U \to\R$ is called running cost.

It is well-known in the literature that the solution $\bx$ to the optimal control problem~\eqref{eqt:optctrl_unconstrained} satisfies the following Hamiltonian ODE system
\begin{equation} \label{eqt:Hamiltonian_ODE}
\begin{dcases}
\dot{\bx}(s) = \nabla_{\bp} H(\bx(s), \bp(s)), & s\in[0,T],\\
\dot{\bp}(s) = -\nabla_{\bx} H(\bx(s), \bp(s)), & s\in[0,T],
\end{dcases}
\end{equation}
where the function $H\colon \R^{2n} \to \R$ is called Hamiltonian defined by
\begin{equation}\label{eqt:defH_unconstrained}
    H(\bx,\bp) = \sup_{\bv\in U} \{ \langle \bv, \bp\rangle - \ell(\bx,\bv)\}, \quad \forall \bx,\bp \in\Rn.
\end{equation}
A method to solve the optimal control problem~\eqref{eqt:optctrl_unconstrained} in high dimensions is to solve the corresponding Hamiltonian ODE system~\eqref{eqt:Hamiltonian_ODE} instead. However, when the Hamiltonian is too complicated, the Hamiltonian ODEs may be hard to solve. To tackle this difficulty, we propose a method called SympOCnet, which uses the SympNet architecture to solve the optimal control problem~\eqref{eqt:optctrl_unconstrained} in high dimensions. The key idea is to perform a change of variables in the phase space using a symplectic map represented by a SympNet, and then solve the Hamiltonian ODEs in the new coordinate system.

Given a symplectic map $\phi\colon \R^{2n}\to \R^{2n}$, we define the change of variables $(\by,\bq) = \phi(\bx,\bp)$ for any $\bx,\bp\in\Rn$. If the phase trajectory $t\mapsto (\bx(t),\bp(t))$ is the solution to the Hamiltonian ODEs~\eqref{eqt:Hamiltonian_ODE}, then by Theorem~\ref{thm:bkgd_symplectic_HODE}, the phase trajectory $t\mapsto (\by(t),\bq(t))$ under the new coordinate system solves the Hamiltonian ODEs with a new Hamiltonian $\tilde{H}\colon \R^{2n}\to\R$ defined by
\begin{equation*}
    \tilde{H}(\by,\bq) = H(\phi^{-1}(\by,\bq)) \quad \forall \by,\bq\in\Rn.
\end{equation*}
In other words, the function $t\mapsto(\by(t), \bq(t))$ satisfies
\begin{equation}\label{eqt:HamiltonianODE_yq_general}
\begin{dcases}
\dot{\by}(s) = \nabla_{\bq} \tilde{H}(\by(s), \bq(s)), & s\in[0,T],\\
\dot{\bq}(s) = -\nabla_{\by} \tilde{H}(\by(s), \bq(s)), & s\in[0,T].
\end{dcases}
\end{equation}
Here, we assume that the new Hamiltonian $\tilde{H}$ has a simpler form such that the corresponding Hamiltonian ODE system~\eqref{eqt:HamiltonianODE_yq_general} is easier to solve. To be specific, we assume that $\tilde{H}$ does not depend on the variable $\by$, and hence the Hamiltonian ODEs become
\begin{equation}\label{eqt:HamiltonianODE_yq}
    \begin{dcases}
    \dot{\by}(s) = \nabla \tilde{H}(\bq(s)),  & s\in[0,T],\\ 
    \dot{\bq}(s) = 0,   & s\in[0,T].
    \end{dcases}
\end{equation}
The solution $s\mapsto \bq(s)$ is a constant function, and the solution $s\mapsto \by(s)$ is an affine function. Therefore, the solution $s\mapsto (\by(s), \bq(s))$ can be represented using three parameters in $\Rn$, which are denoted by $\bq_0 = \bq(0)$, $\by_0 = \by(0)$ and $\bu = \dot{\by}(s)$. With these three parameters, the solution to~\eqref{eqt:HamiltonianODE_yq} is given by $\by(s) = \by_0 + s\bu$ and $\bq(s) = \bq_0$ for any $s\in [0,T]$.

Under this framework, we explain our proposed method in the remainder of this section. In Section~\ref{sec:SympNet}, we apply the SympOCnet method to approximate the unknown function $\phi$ and the unknown parameters $\bu, \by_0$ and $\bq_0$.
Then, after the training process, to further enhance the accuracy we can also apply the pseudospectral method for post-processing (for low-dimensional problems), as we explain in Section~\ref{sec:PSmethod}.

\subsection{SympOCnet method} \label{sec:SympNet}
In this section, we describe how to apply SympOCnet method to solve the optimal control problem~\eqref{eqt:optctrl_unconstrained}. SympOCnet uses 
a SympNet architecture to represent the inverse of the unknown function $\phi$, which is also a symplectic map. In other words, we approximate $\phi^{-1}$ using a SympNet denoted by $\varphi = (\varphi_1, \varphi_2)$.
Here, $\varphi_1$ denotes the first $n$ components in $\varphi$, which is a map from $(\by,\bq)$ to $\bx$, and $\varphi_2$ denotes the last $n$ components in $\varphi$, which is a map from $(\by,\bq)$ to $\bp$. Other learnable parameters include $\bu, \by_0, \bq_0\in\Rn$ mentioned above, which respectively correspond to $\dot{\by}(s)$, $\by(0)$ and $\bq(s)$ in the new Hamiltonian ODE system~\eqref{eqt:HamiltonianODE_yq}.
With the SympNet $\varphi$ and the variables $\bu, \by_0, \bq_0$, the solution to the original Hamiltonian ODEs~\eqref{eqt:Hamiltonian_ODE} is given by $s\mapsto\varphi(\by_0+s\bu, \bq_0)$.

\begin{figure}[htbp]
   \centering \includegraphics[width=0.8\textwidth]{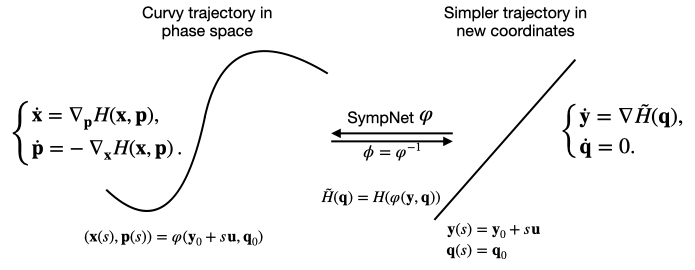}
    \caption{\textbf{An illustration of the SympOCnet method.} We propose to use SympNet $\phi$ to map curvy trajectory in the phase space to simpler trajectory in new coordinates and solve the corresponding Hamiltonian system of the optimal control problem.}
    \label{fig:illustration_sympocnet}
\end{figure}

To learn the SympNet $\varphi$ and the parameters $\bu, \by_0$ and $\bq_0$, we apply the PINN loss~\cite{raissi2019physics} on the Hamiltonian ODEs~\eqref{eqt:Hamiltonian_ODE}, and get
\begin{equation*}
\begin{split}
    \mathcal{L}_{res} 
    &=\sum_{j=1}^N\norm{\frac{d\bx(s_j)}{ds} - \nabla_{\bp} H(\bx(s_j), \bp(s_j)) }^2 + \sum_{j=1}^N\norm{\frac{d\bp(s_j)}{ds} + \nabla_{\bx}H(\bx(s_j), \bp(s_j)) }^2 \\
    &= \sum_{j=1}^N\norm{\nabla_{\by}\varphi_1(\by_0+s_j\bu, \bq_0)\bu - \nabla_{\bp} H(\varphi(\by_0+s_j\bu, \bq_0)) }^2 \\
    &\quad \quad + \sum_{j=1}^N\norm{\nabla_{\by}\varphi_2(\by_0+s_j\bu, \bq_0)\bu + \nabla_{\bx} H(\varphi(\by_0+s_j\bu, \bq_0)) }^2,
\end{split}    
\end{equation*}
where the second equality holds by the following chain rule
\begin{equation*}
\begin{split}
    \frac{d\bx(s_j)}{ds} &= \nabla_{\by}\bx(\by_0+s_j\bu, \bq_0)\frac{d\by(s_j)}{ds} + \nabla_{\bq}\bx(\by_0+s_j\bu, \bq_0)\frac{d\bq(s_j)}{ds}\\ 
    &= \nabla_{\by}\varphi_1(\by_0+s_j\bu, \bq_0)\bu,\\
    \frac{d\bp(s_j)}{ds} &= \nabla_{\by}\bp(\by_0+s_j\bu, \bq_0)\frac{d\by(s_j)}{ds} + \nabla_{\bq}\bp(\by_0+s_j\bu, \bq_0)\frac{d\bq(s_j)}{ds}\\ 
    &= \nabla_{\by}\varphi_2(\by_0+s_j\bu, \bq_0)\bu.
\end{split}
\end{equation*}
Here, the data points are given by $s_1, \dots, s_N \in [0,T]$.
The gradients $\nabla_{\by}\varphi_1(\by_0+s_j\bu, \bq_0)$ and $\nabla_{\by}\varphi_2(\by_0+s_j\bu, \bq_0)$ are calculated by the automatic differentiation of the SympNet $\varphi$.
Moreover, we have another loss term for the initial and terminal conditions of $\bx$, defined by
\begin{equation*}
\begin{split}
    \mathcal{L}_{bd} 
    &= \norm{\bx(0) - \bx_0}^2 +  \norm{\bx(T) - \bx_T}^2\\
    &= \norm{\varphi_1(\by_0,\bq_0) - \bx_0}^2 +  \norm{\varphi_1(\by_0+T\bu,\bq_0) - \bx_T}^2.
\end{split}
\end{equation*}
The total loss function $\mathcal{L}$ is a weighted sum of the two loss terms, defined by
\begin{equation}\label{eqt:def_loss_L}
    \mathcal{L} = \mathcal{L}_{res} + \lambda \mathcal{L}_{bd},
\end{equation}
where $\lambda$ is a tunable positive hyperparameter.
After training, the optimal trajectory is obtained from
\begin{equation}\label{eqt:optimal_x_sympnet}
    \bx(s) = \varphi_1(\by_0 + s\bu, \bq_0) \quad \forall s\in[0,T],
\end{equation}
where the function $\varphi_1$ contains the first $n$ components of the SympNet $\varphi$.
An illustration of our proposed SympOCnet method is shown in Figure~\ref{fig:illustration_sympocnet}.

\subsection{Post-processing using the pseudospectral method}\label{sec:PSmethod}
In Section~\ref{sec:SympNet}, we introduced the SympOCnet to solve the optimal control problem~\eqref{eqt:optctrl_unconstrained}. In this section, we explain how to apply the pseudospectral method to postprocess SympOCnet results to improve feasibility and optimality. 
The pseudospectral method is a popular method proposed to solve a general optimal control problem~\cite{fahroo2001costate, fahroo2002direct, Rao2009Survey,Ross2012Review,garg2011advances}. 
It applies the quadrature rule to discretize the integration and the ODE in a continuous optimal control problem. In this way, the continuous optimal control problem is converted to a finite-dimensional optimization problem with constraints. Then, some optimization algorithms such as sequential quadratic programming are applied to solve the finite-dimensional problem. 

When the dimension is high, it is, in general, hard to solve the converted finite-dimensional optimization problem. The difficulty comes from the number of variables and constraints in the optimization problem, as well as the non-convexity of the problem. However, there are some theoretical guarantees when the initialization is good enough. Therefore, it can be applied as a postprocessing procedure by taking the trajectory obtained from Section~\ref{sec:SympNet} as an initialization.  

We will apply this procedure in the numerical experiments in Sections~\ref{sec:test1_comparison} and~\ref{sec:test2_uncertain}. From the numerical results, we observe that the pseudospectral method preserves the shape of the trajectories obtained from the training step of the SympOCnet method, and it makes the trajectories smoother while improving the optimality. Therefore, the pseudospectral method performs well as a post-processing method. It is also observed in the experiments that pseudospectral method itself with a linear initialization may not converge to the correct solution, which justifies the choice of SympOCnet as a good initializer. Note that we do not apply this post-process to the experiments in Sections~\ref{sec:test3_multidrones} and~\ref{sec:test4_swarm} due to the high dimensionality.
\section{SympOCnet for optimal control problems with state constraints}\label{sec:optctrl_constrained}
In this section, we focus on the following optimal control problem with state constraints
\begin{equation} \label{eqt:optctrl_constrained}
\begin{split}
 &\min\Bigg\{\int_0^T \ell(\bx(s),\bv(s)) ds \colon  h(\bx(s))\geq 0,\dot{\bx}(s)=\bv(s)\in U\forall s\in [0,T],\\
 &\quad\quad\quad\quad\quad\quad\quad\quad\quad\quad\quad\quad\quad\quad\quad\quad\quad\quad\quad\quad  \bx(0)=\bx_0, \bx(T)=\bx_T \Bigg\}\\
    =\,& 
\min\Bigg\{\int_0^T \ell(\bx(s),\dot{\bx}(s)) ds \colon  h(\bx(s))\geq 0, \dot{\bx}(s)\in U\forall s\in [0,T], \\
&\quad\quad \quad\quad\quad\quad\quad\quad\quad\quad\quad\quad\quad\quad\quad\quad\quad
\bx(0)=\bx_0, \bx(T)=\bx_T
\Bigg\},
\end{split}
\end{equation}
where $\ell\colon \Rn\times U \to\R$ is the running cost, and $h\colon \Rn\to \R^m$ is a function providing the constraint on the state variable $\bx$. For instance, to avoid obstacles, $h$ can be defined using the signed distance function to the obstacles.

In order to apply the SympOCnet method in Section~\ref{sec:optctrl_unconstrained}, we need to convert the original problem~\eqref{eqt:optctrl_constrained} to an unconstrained optimal control problem, which is achieved using the soft penalty method described in Section~\ref{subsec:penalty}. 
The soft penalty method will enforce the constraints in an asymptotic sense, i.e., when the hyperparameters converge to zero under some conditions, which is impractical in experiments. In the training process, when the hyperparameters in the soft penalty are fixed, to improve the feasibility of the output trajectory, we introduce extra terms in the loss function in Sections~\ref{sec:training_penalty} and~\ref{sec:training_Lag}. 
The selection of these extra loss terms may differ from problem to problem. Later in Section~\ref{sec:test1_comparison}, we compare the performance of different loss functions under our path planning problem setup to choose a suitable loss function for the other experiments.

\subsection{Soft penalty method for the optimal control problem} \label{subsec:penalty}
To convert a constrained problem to an unconstrained one, the soft penalty method replaces the hard constraint by a soft penalty term in the objective function. 
Here, we consider the log penalty function~\cite{Rucco2015Efficient,Spedicato2018Minimum}, denoted by $\beta_{a}\colon \R\to\R$ and defined by 
\begin{equation}\label{eqt:def_beta_1d}
   \beta_a(x) :=  \begin{dcases}
   -\log(x), & \text{if}\ x>a, \\
   -\log(a) + \frac{1}{2}\left(\left(\frac{x-2a}{a}\right)^2 - 1\right), & \text{if}\ x\leq a,
   \end{dcases}
\end{equation}
where $a$ is a positive hyperparameter.
The high-dimensional log penalty function is the summation of the one-dimensional function acting on each component, i.e., for any $m>1$, the log penalty function $\beta_{a}\colon \R^m\to\R$ is defined by
\begin{equation}\label{eqt:def_beta_hd}
   \beta_a(\bx) :=  \sum_{i=1}^m \beta_a(x_i) \quad \forall \bx=(x_1,\dots,x_m)\in\R^m.
\end{equation}
With this log penalty term, the problem~\eqref{eqt:optctrl_constrained} is converted to
\begin{equation} \label{eqt:optctrl_penalty}
\begin{split}
     \min\Bigg\{\int_0^T
    \ell(\bx(s),\dot{\bx}(s)) + \epsilon \beta_{a}(h(\bx(s))ds \colon \dot{\bx}(s)\in U\forall s\in [0,T],\quad\quad \\
    \bx(0)=\bx_0, \bx(T)=\bx_T\Bigg\},   
\end{split}
\end{equation}
where $\epsilon$ and $a$ are positive hyperparameters. This problem~\eqref{eqt:optctrl_penalty} is in the form of~\eqref{eqt:optctrl_unconstrained}, and hence can be solved using the SympOCnet method in Section~\ref{sec:optctrl_unconstrained}.
By straightforward calculation, for a sequence of parameters $\epsilon_k$ and $a_k$ such that they both converge to zero and the following equations hold
\begin{equation*}
 \lim_{k\to\infty} \epsilon_k\log a_k = 0,\quad\quad 
 \lim_{k\to\infty} \frac{a_k^2}{\epsilon_k} = 0,
\end{equation*}
the penalty term $\epsilon_k\beta_{a_k}$ converges to the indicator function of $[0,+\infty)$, which takes the value $0$ in $[0,+\infty)$ and $+\infty$ in $(-\infty,0)$.
Therefore, the objective function in~\eqref{eqt:optctrl_penalty} with parameters $\epsilon_k$, $a_k$ converges to the objective function in~\eqref{eqt:optctrl_constrained} with the hard constraint given by $h$.

With this penalty term, the Hamiltonian corresponding to the optimal control problem~\eqref{eqt:optctrl_penalty} equals
\begin{equation}\label{eqt:defH_penalty}
H_{\epsilon,a}(\bx,\bp) = \sup_{\bv\in U} \{ \langle \bv, \bp\rangle - \ell(\bx,\bv) - \epsilon \beta_{a}(h(\bx))\}
= H(\bx,\bp) - \epsilon \beta_{a}(h(\bx)),
\end{equation}
where $H$ is the Hamiltonian corresponding to the running cost $\ell$ defined in~\eqref{eqt:optctrl_penalty}.

\begin{remark}
Note that another widely-used method to covert a constrained optimization problem to an unconstrained one is the augmented Lagrangian method. In our problem, since the constraint is enforced for any time $s\in [0,T]$, to apply the augmented Lagrangian method, the Lagrange multiplier needs to be a function of time. Then, with the added Lagrangian term, the new running cost also depends on the time variable, which requires a more complicated symplectic method to handle time-dependent Hamiltonians. Therefore, we do not apply the augmented Lagrangian method for the constrained optimal control problem~\eqref{eqt:optctrl_constrained} in this paper. Instead, we add extra terms in the loss function (see Sections~\ref{sec:training_penalty} and~\ref{sec:training_Lag}) to enforce the constraint. 
\end{remark}

\subsection{Penalty method in the training process of SympOCnet}\label{sec:training_penalty}
To enforce the state constraint $h(\bx(s))\geq 0$, one way is to add the penalty term in the loss function of the training process. 
Here, we introduce two penalty terms, which include the log penalty $\beta_a$ defined in~\eqref{eqt:def_beta_hd} and the quadratic penalty. 
The loss function with the log penalty term corresponding to the state constraint $h(\bx(s))\geq 0$ is defined by
\begin{equation}\label{eqt:def_L_log}
\begin{split}
    \mathcal{L}_{log}(\varphi, \by_0, \bq_0, \bu) &= \mathcal{L}(\varphi, \by_0, \bq_0, \bu) + \frac{\tilde{\lambda}}{N}\sum_{j=1}^N \epsilon \beta_a(h(\varphi_1(\by_0 + s_j \bu, \bq_0)))  \\
     &\quad\quad + \epsilon \beta_a(\bx_0 - \varphi_1(\by_0, \bq_0)) + \epsilon \beta_a(\varphi_1(\by_0, \bq_0) - \bx_0)\\
    &\quad\quad + \epsilon \beta_a(\bx_T - \varphi_1(\by_0 + T \bu, \bq_0)) + \epsilon \beta_a(\varphi_1(\by_0 + T \bu, \bq_0) - \bx_T),
\end{split}
\end{equation}
where $\beta_a$ is the penalty function defined in~\eqref{eqt:def_beta_hd}, $\mathcal{L}$ is the loss function defined in Section~\ref{sec:SympNet}, and $\tilde{\lambda}$ is a positive hyperparameter.
Similarly, the loss function with a quadratic penalty is defined by
\begin{equation}\label{eqt:def_L_quad}
    \mathcal{L}_{quad}(\varphi, \by_0, \bq_0, \bu) = \mathcal{L}(\varphi, \by_0, \bq_0, \bu)
    + \frac{\tilde{\lambda}}{N} \sum_{j=1}^N \|\min\{h(\varphi_1(\by_0 + s_j \bu, \bq_0)), 0\}\|^2,
\end{equation}
where we do not add the penalty term for the initial and terminal conditions, since the quadratic penalty for them are the same with $\mathcal{L}_{bd}$ in the loss function $\mathcal{L}$.

\subsection{Augmented Lagrangian method in the training process of SympOCnet}\label{sec:training_Lag}
Another method to enforce the state constraint $h(\bx(s))\geq 0$ is to apply augmented Lagrangian method in the training process.
The loss function with the Lagrangian term is defined as follows
\begin{equation}\label{eqt:def_L_lag}
\begin{split}
    &\mathcal{L}_{aug}(\varphi, \by_0,\bq_0,\bu, \bmu, \blambda_1, \blambda_2) \\
    =\,& \mathcal{L}(\varphi, \by_0,\bq_0,\bu) + \frac{1}{2\rho_1} \sum_{j=1}^N\|\max\{\bzero, \bmu(s_j) - \rho_1 h(\varphi_1(\by_0 + s_j\bu, \bq_0))\}\|^2\\
    &\quad + \frac{1}{2\rho_2} \left(\|\blambda_1 - \rho_2(\bx_0 - \varphi_1(\by_0, \bq_0))\|^2 + \|\blambda_2 - \rho_2(\bx_T - \varphi_1(\by_0+T\bu, \bq_0))\|^2\right),
\end{split}
\end{equation}
where $\bmu(s_j)\in\R^m$, $\blambda_1,\blambda_2\in \R^n$ are the augmented Lagrange multipliers for the constraints $h(\bx(s_j))\geq 0$, $\bx(0)=\bx_0$ and $\bx(T)=\bx_T$, respectively, and $\rho_1,\rho_2>0$ are positive hyperparameters in the augmented Lagrangian method. The training process becomes an iterative scheme where the $k$-th iteration contains the following two steps
\begin{itemize}
\item Train the SympNet $\varphi$ and parameters $\by_0,\bq_0,\bu\in\Rn$ for several iterations using the loss function $\mathcal{L}_{aug}(\varphi, \by_0,\bq_0,\bu, \bmu^k, \blambda_1^k, \blambda_2^k)$;
\item Update the Lagrange multiplier by 
\begin{equation*}
\begin{split}
 \bmu^{k+1}(s_j) &= \max\{\bzero, \bmu^k(s_j) - \rho_1 h(\varphi_1(\by_0+s_j\bu, \bq_0))\}\quad\forall j=1,\dots,N,\\
 \blambda_1^{k+1} &= \blambda_1^k - \rho_2(\bx_0 - \varphi_1(\by_0, \bq_0)),\\
 \blambda_2^{k+1} &= \blambda_2^k - \rho_2(\bx_T - \varphi_1(\by_0 + T\bu, \bq_0)).
\end{split}
\end{equation*}
\end{itemize}

Note that in practice, we need to tune the hyperparameters $\rho_1$ and $\rho_2$. In our implementation, these hyperparameters are tuned automatically using a similar strategy as in~\cite{Conn1991Globally}.

\section{Applications in path planning problems with obstacle and collision avoidance}\label{sec:numerical_results}
In this section, we apply our method to path planning problems with multiple drones. 
We assume that each drone is represented by a ball in the physical space $\R^m$ with radius $C_{d}$. Note that the meaning of the notation $m$ in this section (i.e., the dimension of the physical space) is different from $m$ in Section~\ref{sec:optctrl_constrained} (which is the number of state constraints).
We set the state variable to be $\bx = (\bx_1,\dots,\bx_M)\in\R^{Mm}$, where $M$ is the number of drones, and each $\bx_j\in\R^m$ denotes the position of the center of each drone.
In other words, the state variable is the concatenation of the position variables of all drones.
Similarly, the control variable $\bv$ is the concatenation of the velocity variables for all drones, i.e., $\bv$ equals $(\bv_1,\dots, \bv_M)\in U^M$, where each $\bv_j$ is the velocity of the $j$-th drone, and $U\subseteq \R^m$ is the space for all admissible velocities for each drone. In our numerical experiments, the control space $U$ is the ball in $\R^m$ with zero center and radius $C_{\bv}>0$, i.e., the norm of the velocity of each drone has the upper bound $C_{\bv}$.

We want to solve for the optimal trajectory with minimal energy, and hence the running cost is set to be
\begin{equation*}
    \ell(\bx,\bv) = \frac{1}{2}\|\bv\|^2\quad \forall \bx\in\R^{Mm}, \bv\in U^M.
\end{equation*}
The corresponding Hamiltonian equals
\begin{equation*}
\begin{split}
    H(\bx,\bp) &= \sup_{\bv\in U^M} \left\{\langle \bv, \bp\rangle - \frac{1}{2}\|\bv\|^2\right\} = \sum_{i=1}^M\sup_{\bv_i\in U} \left\{\langle \bv_i, \bp_i\rangle - \frac{1}{2}\|\bv_i\|^2\right\} \\
    &= \sum_{i=1}^M\begin{dcases}
    \frac{1}{2}\|\bp_i\|^2 & \text{if }\|\bp_i\| \leq C_{\bv},\\
    C_{\bv}\|\bp_i\| - \frac{C_{\bv}^2}{2} & \text{if } \|\bp_i\| > C_{\bv},
    \end{dcases}
\end{split}
\end{equation*}
for any $\bp=(\bp_1,\cdots,\bp_M)\in \R^{Mm}$, where each $\bp_i$ is the momentum variable in $\R^m$ for the $i$-th drone.

To avoid obstacles and collisions among drones, we set the constraint function $h$ to be $h=(h_1,h_2)$, where $h_1$ is for obstacle avoidance, and $h_2$ is for avoiding collisions among drones. If there are $n_o$ obstacles, denoted by $E_1,\dots, E_{n_o}$, then the function $h_1\colon \R^{Mm}\to\R^{Mn_o}$ is defined by
\begin{equation}\label{eqt:def_h1}
    h_1(\bx_1, \dots, \bx_M) = (d(\bx_1), d(\bx_2), \dots, d(\bx_M)) \quad\forall \bx_1,\dots,\bx_M\in\R^m,
\end{equation}
where the function $d\colon \R^m\to\R^{n_o}$ is defined by
\begin{equation}\label{eqt:def_d}
    d(\bx) = (d_1(\bx), \dots, d_{n_o}(\bx))\quad\forall \bx\in\R^m,
\end{equation}
and each function $d_j$ is defined such that $d_j(\bx)<0$ implies the collision of the drone whose center is at the position $\bx$ with the $j$-th obstacle $E_j$.
For instance, $d_j$ can be defined using the signed distance function to $E_j$. The definition of the function $d_j$ depends on the shape of the constraint set $E_j$, and hence we do not specify it now but postpone its definition to each example in later sections.

The function $h_2\colon \R^{Mm}\to\R^{M(M-1)/2}$ is the constraint function for collision avoidance; each component of the constraint function gives a constraint for avoiding collision between a pair of drones. 
Since each drone can be seen as a ball centered at a point in $\R^m$, two drones collide if and only if the distance between their centers is less than the sum of their radii.
Hence, we set the $k$-th component of $h_2$ to be
\begin{equation}\label{eqt:def_h2_distsq}
(h_2)_k(\bx_1,\dots,\bx_M) = \|\bx_i-\bx_j\|^2 - (2C_{d})^2 \quad\forall \bx_1,\dots,\bx_M\in\R^m,
\end{equation}
where $C_d$ is the radius of a drone, and $k=i+(j-1)(j-2)/2$ for any $1\leq i < j \leq M$ is the corresponding constraint index for the pair $(i,j)$. Note that the constraints may be duplicated to simplify the implementation.
As a result, $(h_2)_k(\bx_1,\dots,\bx_M)< 0$ is equivalent to the collision between the $i$-th and $j$-th drones. In other words, such defined constraint function $h_2$ can avoid collisions among drones. 

In the following sections, we present several numerical results. In Section~\ref{sec:test1_comparison}, we compare the performance of different loss functions $\mathcal{L}$, $\mathcal{L}_{log}$, $\mathcal{L}_{quad}$ and $\mathcal{L}_{aug}$ to choose a suitable loss function for the remaining experiments. In Section~\ref{sec:test2_uncertain}, we present an offline-online training strategy for an optimal control problem whose initial conditions are not known in the training process of the SympNet. This experiment shows the generalizability of SympOCnet.
Then, in Section~\ref{sec:test3_multidrones}, we demonstrate the performance of SympOCNet on a high-dimensional problem reported in~\cite{Robinson2018Efficient}. From the results in Section~\ref{sec:test3_multidrones}, we observe that SympOCnet can handle the path planning problem whose state space has dimension $512$, and hence the method can potentially mitigate the CoD.  In Section~\ref{sec:test4_swarm}, we apply SympOCnet to a swarm path planning problem in~\cite{onken2021neural}, and demonstrate good performance and efficiency in path planning problems, where the agents move in a three-dimensional space.

In Sections~\ref{sec:test1_comparison} and~\ref{sec:test2_uncertain} we use a $6$-layer SympNet with $60$ hidden neurons in each layer and ReLU activation function. In Sections~\ref{sec:test3_multidrones} and~\ref{sec:test4_swarm} we use a $6$-layer SympNet with $200$ hidden neurons in each layer and ReLU activation function. In all these experiments, the parameter $\lambda$ in~\eqref{eqt:def_loss_L} is set to be $600$, and the parameter $\tilde{\lambda}$ in~\eqref{eqt:def_L_log} and~\eqref{eqt:def_L_quad} is set to be $200$. By default, the time horizon $T$ is set to be $1$, and the speed limit $C_{\bv}$ is $25$.
All the numerical experiments are run using a shared NVIDIA GeForce RTX 3090 GPU and a shared NVIDIA RTX A6000 GPU.
In each experiment, we only show the figure for a part of the testing cases. More numerical results and codes are provided in \url{https://github.com/zzhang222/SympOCNet}. Video animations of these examples are available online at \url{https://github.com/zzhang222/SympOCNet/tree/main/saved_results_animation}.

\subsection{Comparison among different loss functions}\label{sec:test1_comparison}

\begin{figure}[htbp]
   \centering \includegraphics[width=0.25\textwidth]{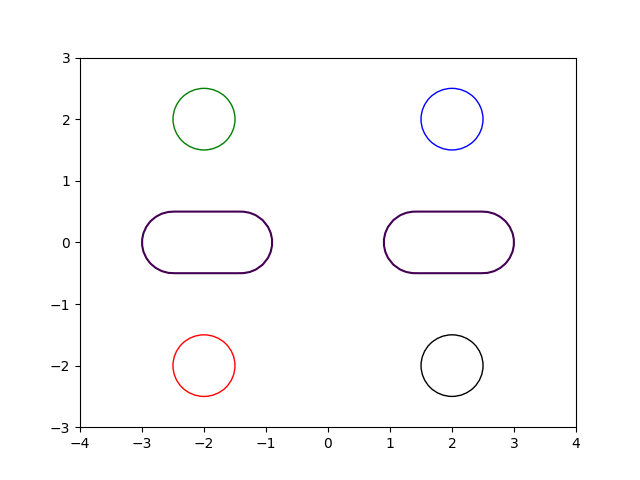}
    \caption{\textbf{Initial positions of drones and obstacles.} The four drones located at their initial positions are represented by the four colored circles. The two obstacles are represented by round cornered rectangles in the middle.}
    \label{fig:test1_initpos}
\end{figure}

In this section, we consider the path planning problem with two obstacles and four drones in $\R^2$. Each obstacle $E_j\subset \R^2$ is a convex set containing points whose distance to a line segment (denoted by $l_j\subset \R^2$) is less than a constant $C_o$. The initial positions of the four drones are $(-2, -2)$, $(2, -2)$, $(2, 2)$, $(-2, 2)$, and their terminal positions are $(2, 2), (-2, 2), (-2, -2), (2, -2)$. Therefore, the optimal control problem~\eqref{eqt:optctrl_constrained} is an $8$-dimensional problem, whose initial state $\bx_0$ is $(-2, -2, 2, -2, 2, 2, -2, 2)$ and terminal state $\bx_T$ is $(2, 2, -2, 2, -2, -2, 2, -2)$. The obstacles and the initial positions of the four drones are shown in Figure~\ref{fig:test1_initpos}. 

\begin{figure}[htbp]
    \centering
    \begin{subfigure}{0.16\textwidth}
        \centering \includegraphics[width=\textwidth]{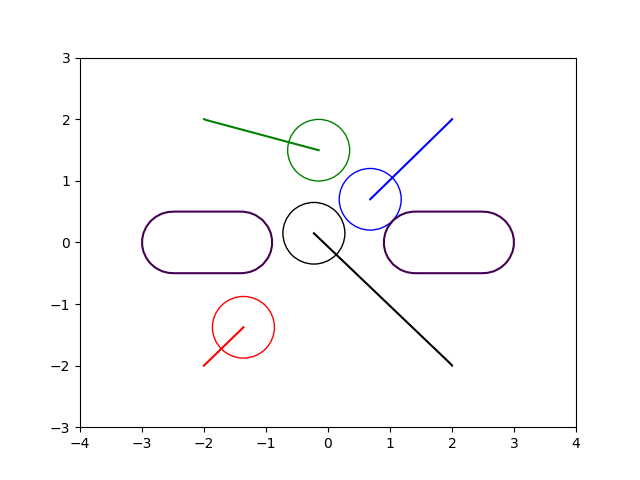}
        \caption{NN, $t=\frac{1}{3}$}
    \end{subfigure}
    \begin{subfigure}{0.16\textwidth}
        \centering \includegraphics[width=\textwidth]{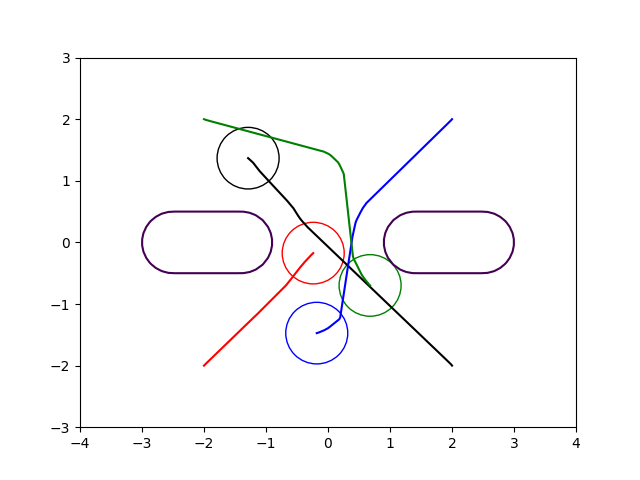}
        \caption{NN, $t=\frac{2}{3}$}
    \end{subfigure}
    \begin{subfigure}{0.16\textwidth}
        \centering \includegraphics[width=\textwidth]{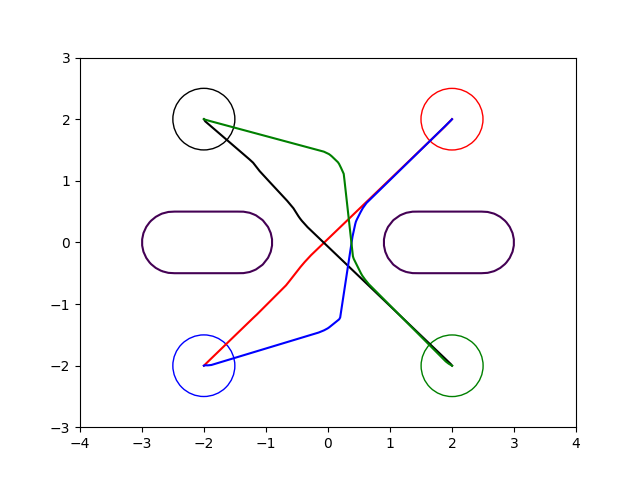}
        \caption{NN, $t=1$}
    \end{subfigure}
    \begin{subfigure}{0.16\textwidth}
        \centering \includegraphics[width=\textwidth]{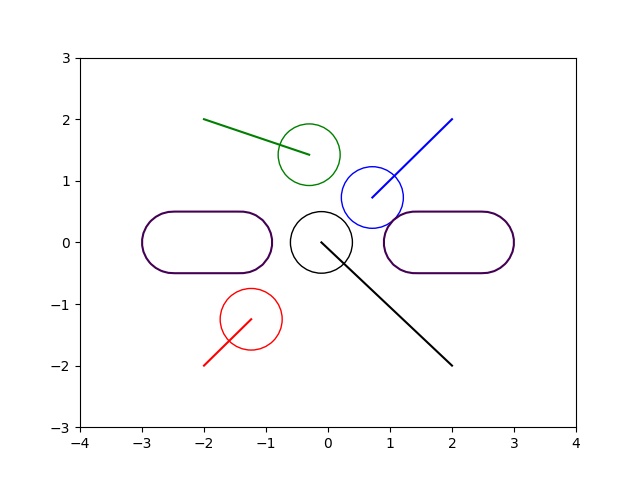}
        \caption{PS, $t=\frac{1}{3}$}
    \end{subfigure}
    \begin{subfigure}{0.16\textwidth}
        \centering \includegraphics[width=\textwidth]{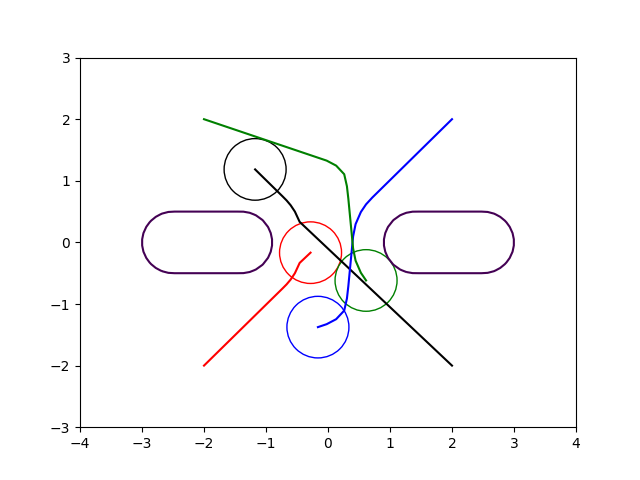}
        \caption{PS, $t=\frac{2}{3}$}
    \end{subfigure}
    \begin{subfigure}{0.16\textwidth}
        \centering \includegraphics[width=\textwidth]{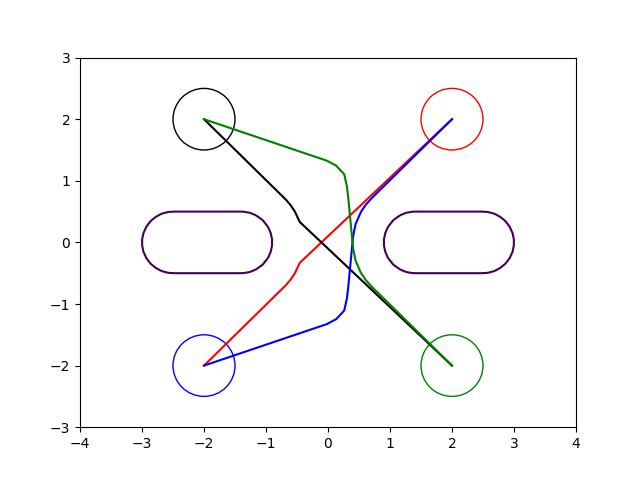}
        \caption{PS, $t=1$}
    \end{subfigure}
    \caption{\textbf{Output trajectories with loss $\mathcal{L}_{aug}$ and pseudospectral method.} The left three figures show the trajectories of our proposed SympOCnet method at different time $t=\frac{1}{3}, \frac{2}{3},1$, and the right three figures show the outputs of the post-process (pseudospectral method) at different time $t=\frac{1}{3}, \frac{2}{3},1$. In each figure, the four colored circles show the positions of the four drones at time $t$, and the curves connected to the circles show the trajectories before time $t$. }
    \label{fig:test1_distsq}
\end{figure}

To avoid the collisions of drones with obstacles, we define the first part of the constraint function $h$ by~\eqref{eqt:def_h1}, where each component $d_j$ of the function $d$ in~\eqref{eqt:def_d} is defined by
\begin{equation*}
    d_j(\bx) = \min_{\by\in l_j}\|\bx-\by\|^2 - (C_o+C_d)^2 \quad \forall \bx\in\R^2.
\end{equation*}
The term $\min_{\by\in l_j}\|\bx-\by\|^2$ in $d_j$ gives the squared distance between $\bx$ and the line segment $l_j$. Here, we use the squared distance instead of distance function to improve the smoothness of the constraint function $h$.
Since the obstacle $E_j$ contains all points whose distance to the line segment $l_j$ is less than $C_o$, a drone collide with $E_j$ if and only if the distance between its center and the line segment $l_j$ is less than $C_o + C_d$. Therefore, such defined function $d_j$ provides a constraint which avoids the collisions of drones with the $j$-th obstacle $E_j$.

Under this problem setup, we compare four different loss functions, i.e., the function $\mathcal{L}$ defined in~\eqref{eqt:def_loss_L}, $\mathcal{L}_{log}$ in~\eqref{eqt:def_L_log}, $\mathcal{L}_{quad}$ in~\eqref{eqt:def_L_quad}, and $\mathcal{L}_{aug}$ in~\eqref{eqt:def_L_lag}. With each loss function, we train the neural network for $100,000$ iterations, and then run the pseudospectral method to improve the results. By comparing the performance of the four loss functions in this example, we choose a better loss function for the other experiments in the remainder of this paper.

The outputs of the four loss functions are shown in Table~\ref{tab:test1_dist_square}, where the hyperparameters in the penalty term $\epsilon \beta_{a}(h(\bx(s))$ are set to be $a = 0.004$ and $\epsilon = 0.0004$. 
Results for different loss functions are shown in different columns. For each loss function, we repeat the experiment $10$ times with different random seeds, and then take average to obtain the statistics in the table. We observe the convergence of the pseudospectral method in all repeated experiments, which shows the robustness of our SympOCnet method, providing optimal or suboptimal solutions in this example.
The second to fourth lines in Table~\ref{tab:test1_dist_square} show the minimal constraint value, the cost value in the optimal control problem, and the running time of the training process of our proposed SympOCnet method. 
The last four lines in Table~\ref{tab:test1_dist_square} show the minimal constraint value, the cost value in the optimal control problem, the number of iterations, and the running time of the post-process using pseudospectral method.

\begin{table}[htbp]
    \centering
    \begin{tabular}{c|c|c|c|c|c}
        \hline
        \multicolumn{2}{c|}{loss function} & $\mathcal{L}$  & $\mathcal{L}_{log}$ & $\mathcal{L}_{quad}$ & $\mathcal{L}_{aug}$ \\
        \hline
        \hline
        \multirow{3}{*}{SympOCnet} & min constraint & -0.834  &  -0.00704 & -0.0824 & -4.84E-04 \\
        & cost & 79.62  & 94.15 & 76.50 & 93.98 \\
        & running time (s) &  1229 & 1776 & 1361 & 1708 \\
        \hline
        \hline
        \multirow{4}{*}{PS} & min constraint & -4.70E-09 & -9.87E-12 & -2.04E-12 & -8.02E-13 \\
        & cost & 73.47  & 91.00 & 76.95 & 80.91 \\
        & \# iterations & 49 & 19 & 17 & 23 \\
        & running time (s) & 104 & 44 & 38 & 47 \\
        \hline
    \end{tabular}
    \hfill 
    \caption{\textbf{The comparison among results of four loss functions with SympOCnet method and pseudospectral post-process.} We show the minimal constraint value $\min_{s}h(\bx(s))$, the cost $\int_{s=0}^T\frac{||\bv(s)||^2}{2}ds$, and the running time of the SympOCnet as well as those statistics after the post-process. The loss function $\mathcal{L}_{aug}$ provides the solution with least amount of constraint violations and reasonable cost values. It can be seen that in all cases, the SympOCnet provides a good initialization for the pseudospectral method. }
    \label{tab:test1_dist_square}
\end{table}

\begin{table}[htbp]
    \centering
    \begin{tabular}{c|c|c|c}
        \hline
        minimal constraint  & cost & number of iterations & running time (s) \\
        \hline
        \hline
        -1.00E+00 & 64.00 & 5 & 44\\
        \hline
    \end{tabular}
    \hfill 
    \caption{\textbf{The result of the pseudospectral method with the linear initialization.} We show the minimal constraint value $\min_{s}h(\bx(s))$, the cost $\int_{s=0}^T\frac{||\bv(s)||^2}{2}ds$, the number of iterations, and the running time of the pseudospectral method with the linear initialization (i.e., the initial trajectory is an affine function of the time variable, and the initial velocity is a constant function). The minimal constraint value is $-1.00$, which shows that the constraint is not satisfied, and the pseudospectral method with the linear initialization does not converge to a feasible solution.}
    \label{tab:test1_ps_linearinit}
\end{table}

Comparing the results of the four loss functions, we observe that the penalty term and the augmented Lagrangian term in the loss function indeed improve the feasibility of the obtained trajectory. The constraints are better satisfied using the loss $\mathcal{L}_{aug}$ in~\eqref{eqt:def_L_lag}, while the corresponding cost value may be a little bit higher than the the results of $\mathcal{L}_{quad}$. In the other experiments, the problems are more complicated, and the constraints are harder to satisfy. Therefore, we use the loss function $\mathcal{L}_{aug}$ for the experiments in the remainder of the paper to help enforce the state constraints. 
For simpler problems, other loss function may be more preferable for faster running time. How to choose the loss function adaptively may be a possible future direction. 
Although the augmented Lagrangian method and log penalty loss function are slightly slower than the other two loss functions, it takes less than $30$ minutes to run $100,000$ training iterations, which shows the efficiency of SympOCnet. 

To show the improvement using our proposed SympOCnet method as an initialization of the pseudospectral method, we solve the same problem using the pseudospectral method with the linear initialization, i.e., we set the initial guess of the pseudospectral method to be $\bv(t) = \frac{\bx_T- \bx_0}{T}$ and $\bx(t) = \frac{(\bx_T- \bx_0)t}{T}$ for any $t\in [0,T]$. The minimal constraint value, the cost value in the optimal control problem, the number of iterations, and the running time are shown in Table~\ref{tab:test1_ps_linearinit}. The minimal constraint value in Table~\ref{tab:test1_ps_linearinit} is $-1.00$, which shows that the pseudospectral method with the linear initialization does not converge to a feasible solution. Therefore, compared with the linear initialization, our SympOCnet method, no matter which loss function in $\mathcal{L}$, $\mathcal{L}_{log}$, $\mathcal{L}_{quad}$, and $\mathcal{L}_{aug}$ we choose, provides a better initial guess for the pseudospectral method. 

Moreover, we show the obtained results in one trial using SympOCnet with loss $\mathcal{L}_{aug}$ and the post-process using the pseudospectral method in Figure~\ref{fig:test1_distsq}. The left three figures show the trajectories of our proposed SympOCnet method at different times, while the right three figures correspond to the final results from the pseudospectral method at different times. In each figure, the four colored circles show the positions of the four drones at the specific time, and the curves connected to the circles show the trajectories before the current time. We observe that the shape of the trajectory given by our SympOCnet method does not change too much after the post-process, which shows that SympOCnet provides a suboptimal solution to the optimal control problem in this example.

\subsection{Generalizability and offline-online training}\label{sec:test2_uncertain}

\begin{figure}[htbp]
    \centering
    \begin{subfigure}{0.16\textwidth}
        \centering \includegraphics[width=\textwidth]{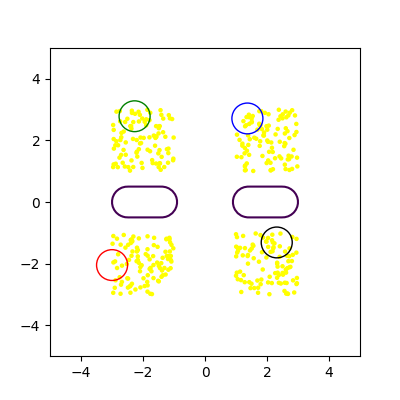}
        \caption{case 1}
    \end{subfigure}
    \begin{subfigure}{0.16\textwidth}
        \centering \includegraphics[width=\textwidth]{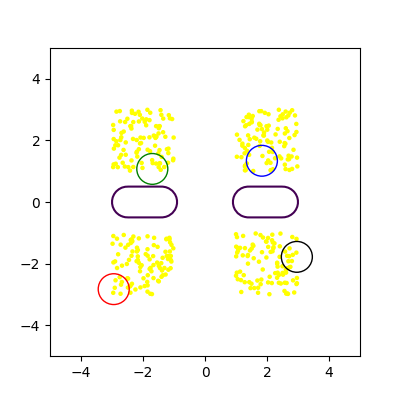}
        \caption{case 2}
    \end{subfigure}
    \begin{subfigure}{0.16\textwidth}
        \centering \includegraphics[width=\textwidth]{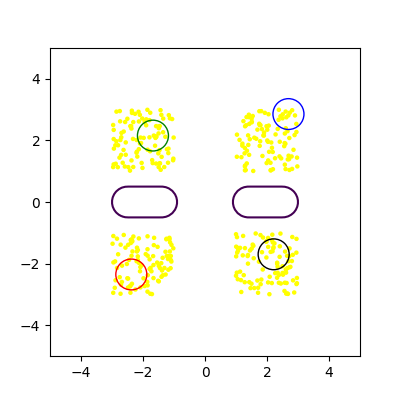}
        \caption{case 3}
    \end{subfigure}
    \begin{subfigure}{0.16\textwidth}
        \centering \includegraphics[width=\textwidth]{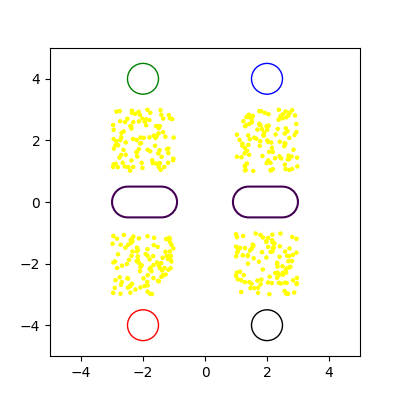}
        \caption{case 4}
    \end{subfigure}
    \begin{subfigure}{0.16\textwidth}
        \centering \includegraphics[width=\textwidth]{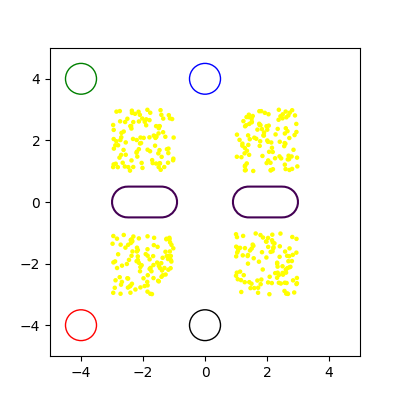}
        \caption{case 5}
    \end{subfigure}
    \begin{subfigure}{0.16\textwidth}
        \centering \includegraphics[width=\textwidth]{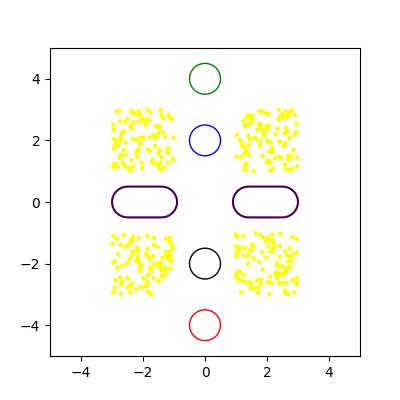}
        \caption{case 6}
    \end{subfigure}
    \caption{\textbf{Initial positions of drones, obstacles, and training data.} In the six figures, we show the initial positions of the four drones (represented by colored circles), the positions of the two obstacles (represented by the round cornered rectangles), and the randomly sampled training data (represented by the yellow dots) in the six testing cases in Section~\ref{sec:test2_uncertain}.}
    \label{fig:test2_IC}
\end{figure}

\begin{figure}[htbp]
    \centering
    \begin{subfigure}{0.16\textwidth}
        \centering \includegraphics[width=\textwidth]{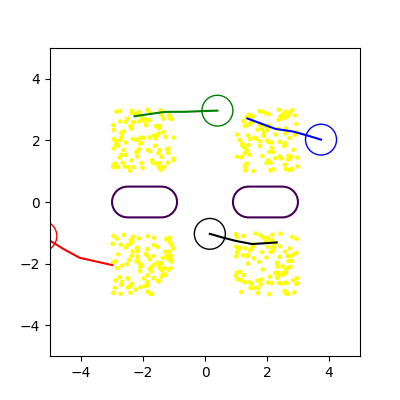}
        \caption{NN, $t=\frac{1}{3}$}
    \end{subfigure}
    \begin{subfigure}{0.16\textwidth}
        \centering \includegraphics[width=\textwidth]{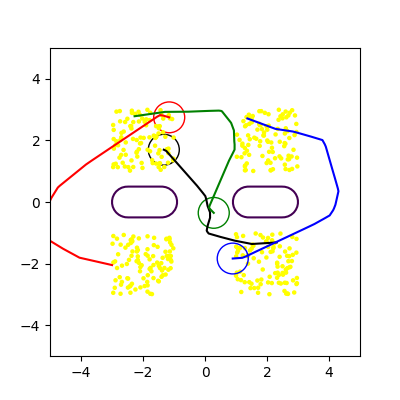}
        \caption{NN, $t=\frac{2}{3}$}
    \end{subfigure}
    \begin{subfigure}{0.16\textwidth}
        \centering \includegraphics[width=\textwidth]{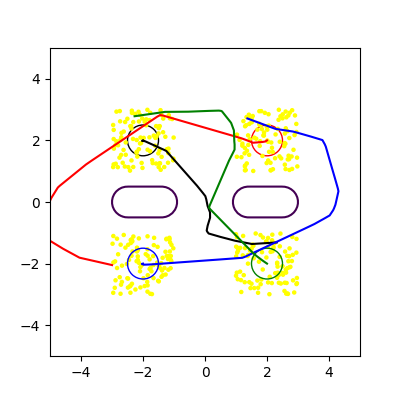}
        \caption{NN, $t=1$}
    \end{subfigure}
    \begin{subfigure}{0.16\textwidth}
        \centering \includegraphics[width=\textwidth]{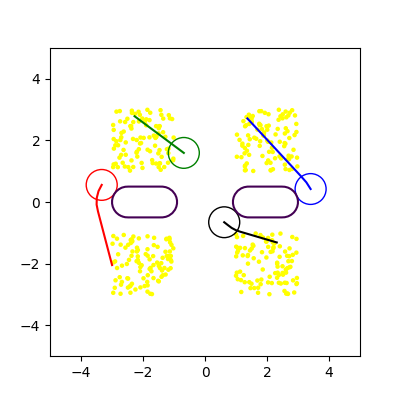}
        \caption{PS, $t=\frac{1}{3}$}
    \end{subfigure}
    \begin{subfigure}{0.16\textwidth}
        \centering \includegraphics[width=\textwidth]{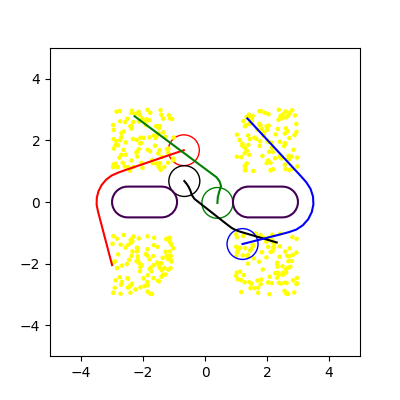}
        \caption{PS, $t=\frac{2}{3}$}
    \end{subfigure}
    \begin{subfigure}{0.16\textwidth}
        \centering \includegraphics[width=\textwidth]{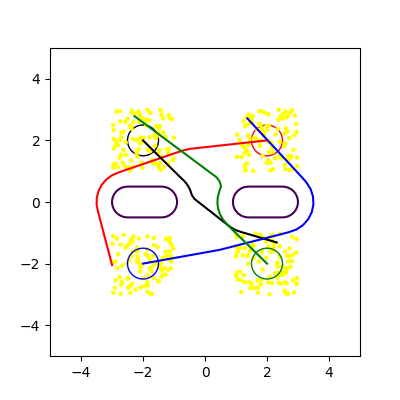}
        \caption{PS, $t=1$}
    \end{subfigure}
    \\
    \begin{subfigure}{0.16\textwidth}
        \centering \includegraphics[width=\textwidth]{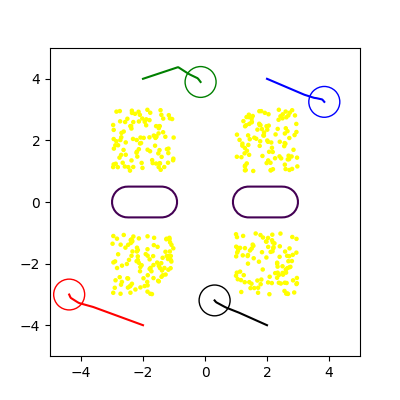}
        \caption{NN, $t=\frac{1}{3}$}
    \end{subfigure}
    \begin{subfigure}{0.16\textwidth}
        \centering \includegraphics[width=\textwidth]{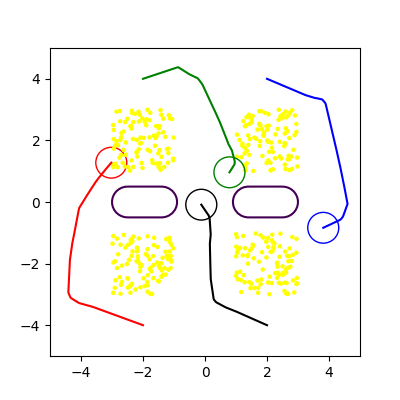}
        \caption{NN, $t=\frac{2}{3}$}
    \end{subfigure}
    \begin{subfigure}{0.16\textwidth}
        \centering \includegraphics[width=\textwidth]{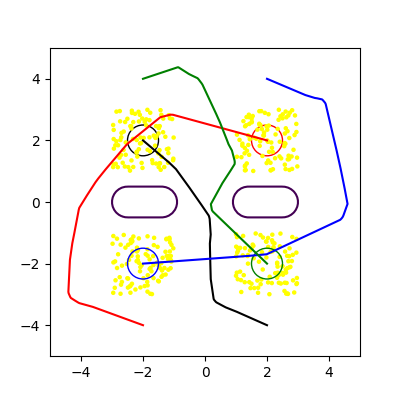}
        \caption{NN, $t=1$}
    \end{subfigure}
    \begin{subfigure}{0.16\textwidth}
        \centering \includegraphics[width=\textwidth]{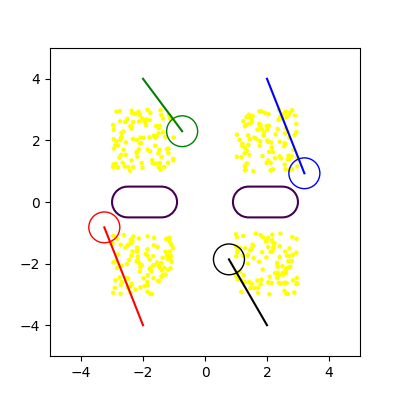}
        \caption{PS, $t=\frac{1}{3}$}
    \end{subfigure}
    \begin{subfigure}{0.16\textwidth}
        \centering \includegraphics[width=\textwidth]{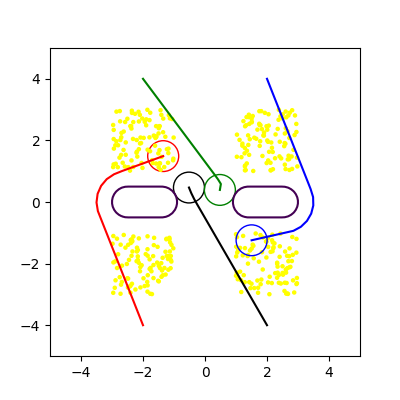}
        \caption{PS, $t=\frac{2}{3}$}
    \end{subfigure}
    \begin{subfigure}{0.16\textwidth}
        \centering \includegraphics[width=\textwidth]{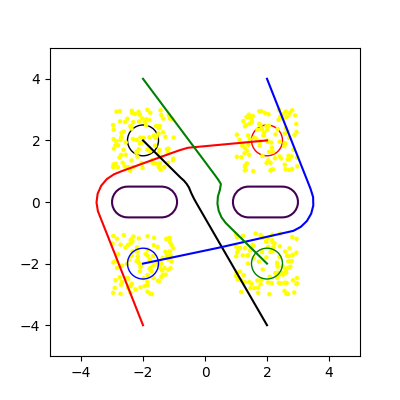}
        \caption{PS, $t=1$}
    \end{subfigure}
    \caption{\textbf{Output trajectories of offline-online training and post-process.} The first and second rows correspond to the first and fourth case, respectively. The left three columns are the output from L-BFGS, and the right three columns are the output from the pseudospectral method. In each figure, the four colored circles show the positions of the four drones at time $t$, and the curves connected to the circles show the trajectories before time $t$. }
    \label{fig:test2}
\end{figure}

In this section, we consider the same problem setup as in Section~\ref{sec:test1_comparison}. To test the generalizability of SympOCnet, we assume that the exact initial position is unknown in the training process. We train one SympNet $\varphi$ using multiple initial positions sampled from a uniform distribution centered at the initial position in Section~\ref{sec:test1_comparison}. With multiple initial positions, we train one SympNet $\varphi$ and several variables $\{(\by_{0k},\bu_k,\bq_{0k})\}_k$, where the $k$-th trainable tuple $(\by_{0k},\bu_k,\bq_{0k})$ corresponds to the $k$-th sampled initial position. We refer to this training process of the SympNet and all trainable variables as ``offline training".
After $100,000$ steps of offline training, the exact initial positions are provided. To obtain an improved solution in reasonable running time using the trained SympNet $\varphi$ as well as the new information, we fix the SympNet $\varphi$ and apply the L-BFGS method to train the variables $(\by_0,\bu,\bq_0)$ with the exact initial and terminal positions. We call this training process using L-BFGS method the ``online training''. The optimal trajectory is then computed using~\eqref{eqt:optimal_x_sympnet}, where $\varphi_1$ comes from the SympNet $\varphi$ trained in the offline training, and the parameters $(\by_0,\bu,\bq_0)$ are outputted from the online training. After the online training, we also apply the pseudospectral method to improve the quality of the results.

To generate training data for the initial positions, we sample $100$ points from the uniform distribution in $[x-1, x+1]\times [y-1,y+1]$ where $(x,y)$ is the initial position of a drone in Section~\ref{sec:test1_comparison}. 
To test the generalizability of our method, we have six cases. In the first three cases, the exact initial positions are randomly sampled from the same distribution used to generate the training data. In other words, we assume the exact data is covered in the area of the training data. In the last three cases, we set the initial conditions for the state variable to be $(-2, -4, 2, -4, 2,4, -2,4)$, $(-4, -4, 0, -4, 0,4, -4,4)$, and $(0, -4, 0, -2, 0,2, 0,4)$, respectively. 
These three initial positions are not in the range of the training data. Therefore, the results with these three initial positions show the generalizability of our SympOCnet method when the test data is not covered by the area of the training data. 
The initial positions of the drones, the positions of the obstacles, and the training data are shown in Figure~\ref{fig:test2_IC}, where the four colored circles denote the drones, the two round cornered rectangles denote the obstacles, and the yellow dots denote the sampled training data.

In our experiments, we found that the performance is improved if we increase the dimension of the problem by adding some latent variables. We increase the dimension of the state space from $8$ to $16$, and define the Hamiltonian on the larger space by
\begin{equation*}
    H(\bx_1,\bx_2,\bp_1,\bp_2) = H_{\epsilon,a}(\bx_1,\bp_1) + \frac{1}{2}\|\bp_2\|^2\quad \forall \bx_1,\bx_2,\bp_1,\bp_2\in \R^8,
\end{equation*}
where $(\bx_1,\bp_1)$ are the variables in the original problem, $(\bx_2,\bp_2)$ are the latent variables we add, and $H_{\epsilon,a}$ is the function defined in~\eqref{eqt:defH_penalty}.
The minimal constraint, cost in the optimal control problem, and running time of the L-BFGS method and pseudospectral method are shown in Table~\ref{tab:test2}. As in Section~\ref{sec:test1_comparison}, we run $10$ repeated experiments and put their average values in the table.
In the online training, we run $1,000$ L-BFGS iterations to train the parameters of all the six cases simultaneously. Note that we can train them all at once, since they share the same SympNet, and the loss function is the summation of the loss functions in all cases.
Then, we run the pseudospectral method on each case until it converges. 
The statistics for the L-BFGS training and the pseudospectral method for the six cases are shown in each line of Table~\ref{tab:test2}. 
The minimal constraint value of L-BFGS method is not as good as the results in Section~\ref{sec:test1_comparison}, which is expected, since the initial positions are not known when the SympNet is trained. However, from the minimal constraint values of the pseudospectral post-process in the six cases, we conclude that our offline-online training still provides a good initialization for the pseudospectral method, which shows a good generalizability of our proposed method.

\begin{table}[htbp]
    \centering
    \begin{tabular}{c|c|c|c}
        \hline
         & minimal constraint & cost & time (s) \\
        \hline
        \hline
        L-BFGS & -1.27E-04 & 231.44 & 158\\
        \hline
        PS (case 1) & -1.97E-11 & 102.00 & 59 \\
        PS (case 2) & -1.74E-11 & 94.59 & 77 
        \\
        PS (case 3) & -3.36E-11 & 90.95 & 79
        \\
        PS (case 4) & -6.00E-09 & 141.22 & 65
        \\
        PS (case 5) & -7.19E-08 & 140.47 & 111
        \\
        PS (case 6) & 1.64E-13 & 112.41 & 67
        \\
        \hline
    \end{tabular}
    \hfill 
    \caption{\textbf{The results given by the offline-online training (L-BFGS) and the post-process (pseudospectral method) in six generalizability tests.} The offline-online training using L-BFGS method provides a good initialization for the pseudospectral method. With this initialization, the pseudospectral method in the post-process improves the feasibility and optimality of the solution in all the six cases.}
    \label{tab:test2}
\end{table}

\begin{table}[htbp]
    \centering
    \begin{tabular}{c|c|c|c}
        \hline
         & minimal constraint & cost & time (s) \\
        \hline
        \hline
        PS (linear initialization, case 1) & -7.71E-09 & 85.19 & 134 \\
        PS (linear initialization, case 2) & 2.95E-07 & 83.98 & 173 
        \\
        PS (linear initialization, case 3) & -1.94E-10 & 79.67 & 178
        \\
        PS (linear initialization, case 4) & -5.55E-16 & 80.44 & 157
        \\
        PS (linear initialization, case 5) & -8.08E-14 & 136.16 & 147
        \\
        PS (linear initialization, case 6) & -1.00 & 60.00 & 28
        \\
        \hline
    \end{tabular}
    \hfill 
    \caption{\textbf{The results given by the pseudospectral method with the linear initialization in six generalizability tests.} Although the pseudospectral method with the linear initialization performs slightly better than our proposed method in the first five cases, it violates the constraint in the last case. Therefore, our proposed method provides a more robust initialization for the pseudospectral method compared with the linear initialization.}
    \label{tab:test2_ps_linearinit}
\end{table}

We also compare our results with the pseudospectral method whose initialization is given by the linear interpolation (as descibed in Section~\ref{sec:test1_comparison}) in the six generalizability tests. The minimal constraint value, cost value in the optimal control problem, and running time of the pseudospectral method with the linear initialization are shown in Table~\ref{tab:test2_ps_linearinit}. Comparing Table~\ref{tab:test2} with Table~\ref{tab:test2_ps_linearinit}, we see that although our proposed method has slightly larger cost values in the first five cases, the pseudospectral method with the linear initialization does not converge to a feasible solution in the last case. Therefore, our proposed method provides a more robust initialization for the pseudospectral method compared with the linear initialization.

Moreover, we plot the output trajectories of one trial of the first and fourth cases computed using our proposed method in Figure~\ref{fig:test2}. The left three columns show the output trajectory given by the offline-online training, and the right three columns are the output trajectory from the post-process. In each figure, the yellow dots are the initial positions of the sampled training data.
Similar as in Section~\ref{sec:test1_comparison}, 
the four colored circles show the current positions of the four drones, and the curves connected to them illustrate the trajectories before the current time.
The first row in Figure~\ref{fig:test2} is for the first case where the exact initial positions are in the area of the sampled training data. The second row is for the fourth case where the exact initial positions are outside the area of the training data. We observe that the offline-online training process still provides feasible trajectories in both cases, based on which the post-process improves the smoothness of the trajectories to provide more optimal solutions. Therefore, our proposed SympOCnet method has the generalizability to handle unseen data which are close to the training data.

\subsection{Multiple drones in a two-dimensional room} \label{sec:test3_multidrones}
In this section, we consider the path planning problem for multiple drones in a room, which is inspired by~\cite{Robinson2018Efficient}. To be specific, we consider $M$ drones moving in a room $[-C_r,C_r]^2$ for some $C_r>0$. In our experiments, we set the constant $C_r = 5$ (i.e., the room has size $10\times 10$). We need to avoid the collisions among the drones as well as the collisions between drones and the walls.
The constraint function $h=(h_1,h_2)$ contains two parts: the function $h_2$ is defined by~\eqref{eqt:def_h2_distsq}, and the function $h_1$ is the obstacle constraint provided by the four walls. Hence, we set $h_1$ in the form of~\eqref{eqt:def_h1}, where the function $d$ is defined by
$
d(\bx) = (x_1 + C_r, C_r-x_1, x_2 + C_r, C_r-x_2)$ for any $\bx=(x_1,x_2)\in\R^2$.
The initial positions of the drones are near the boundary of the room, and the terminal positions are the opposite locations, i.e., we set $\bx_T = -\bx_0$. This is a high-dimensional example (the dimension of the state space is $n=2M$, which ranges from $64$ to $512$ in our experiments), and hence we apply our SympOCnet method without the post-process.

First, we set the drone radius $C_d$ to be $0.3$. We run $10$ repeated experiments (with $100,000$ training iterations in each experiment) with drone numbers $M=32$ and $64$ to test the robustness of our SympOCnet method. The results are shown in Table~\ref{tab:test3}. In each trial, we compute the minimal normalized distance of any pair of drones at any time grid, i.e., we compute $\mathcal{D}$ defined by
\begin{equation*}
    \mathcal{D} = \frac{1}{2C_d}\min_{1\leq i<j\leq M} \min_{1\leq k\leq N} \|\bx_i(t_k)-\bx_j(t_k)\|,
\end{equation*}
where $\bx_i(t_k)$ denotes the center position for the $i$-th drone at time $t_k$.
Then, we compute the average and standard deviation of these minimal normalized distances $\mathcal{D}$ in the $10$ trials, which are shown in the second and third columns in Table~\ref{tab:test3}. We also compute the scaled cost in the optimal control problem, which is defined by
$    \frac{1}{M} \int_0^T \frac{1}{2}\|\dot{\bx}(t)\|^2 dt$,
where $\bx(\cdot)$ is the output trajectory for the state variable $\bx$ in the optimal control problem.
The average of the scaled cost is shown in the fourth column in Table~\ref{tab:test3}. In the last column of Table~\ref{tab:test3}, we show the averaged running time of the training process.
From the results, we observe that our SympOCnet method provides feasible results in reasonable running time for this high-dimensional problem. 
Note that the collision between two drones happens whenever the minimal normalized distance is less than $1$. Therefore, in general the collision does not happen for $32$ drones, but it may happen for $64$ or more drones. To provide a more feasible solution, in the following experiments, we put a safe zone outside each drone by setting $C_d$ to be a little bit larger than the actual drone radius.

\begin{table}[htbp]
    \centering
    \begin{tabular}{c|c|c|c|c}
        \hline
        \# drones & $\mathbb{E}(\mathcal{D})$ & std($\mathcal{D}$) & scaled cost & time (s) \\
        \hline
        \hline
        32 & 1.001 & 0.00390 & 85.84 & 1905\\
        64 & 0.986 & 0.00347 & 100.43 & 2214\\
        \hline
    \end{tabular}
    \hfill 
    \caption{\textbf{The results for the path planning problem with $32$ and $64$ drones in a room.} We show the expected value and standard deviation of the minimal normalized distances. The collision does not happen for 32 drones, but it may happen for 64 or more drones. Therefore, we need to put a safe zone near each drone to avoid collision when the number of drones is large. The computation time increases only by 309 seconds when scaling from 32 to 64 drones. }
    \label{tab:test3}
\end{table}

\begin{figure}[htbp]
    \centering
    \begin{subfigure}{0.24\textwidth}
        \centering \includegraphics[width=\textwidth]{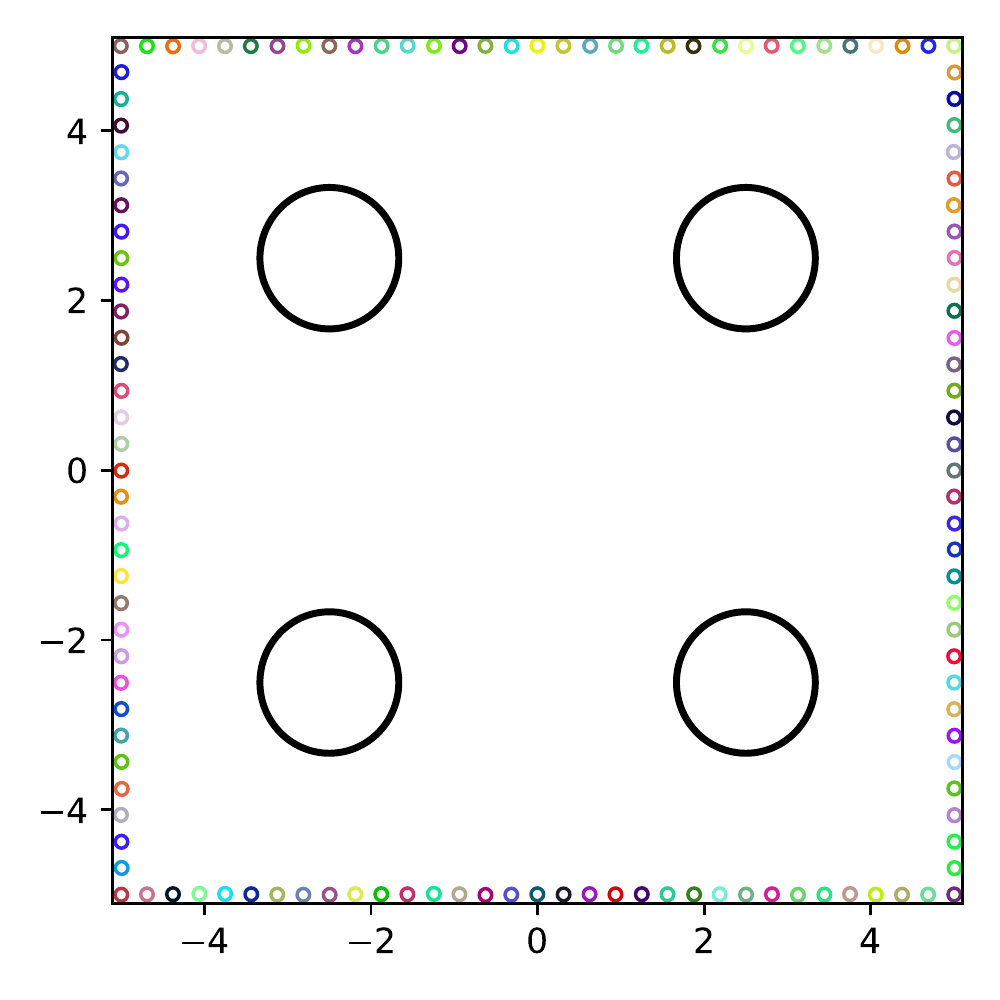}
        \caption{$t=0$}
    \end{subfigure}
    \begin{subfigure}{0.24\textwidth}
        \centering \includegraphics[width=\textwidth]{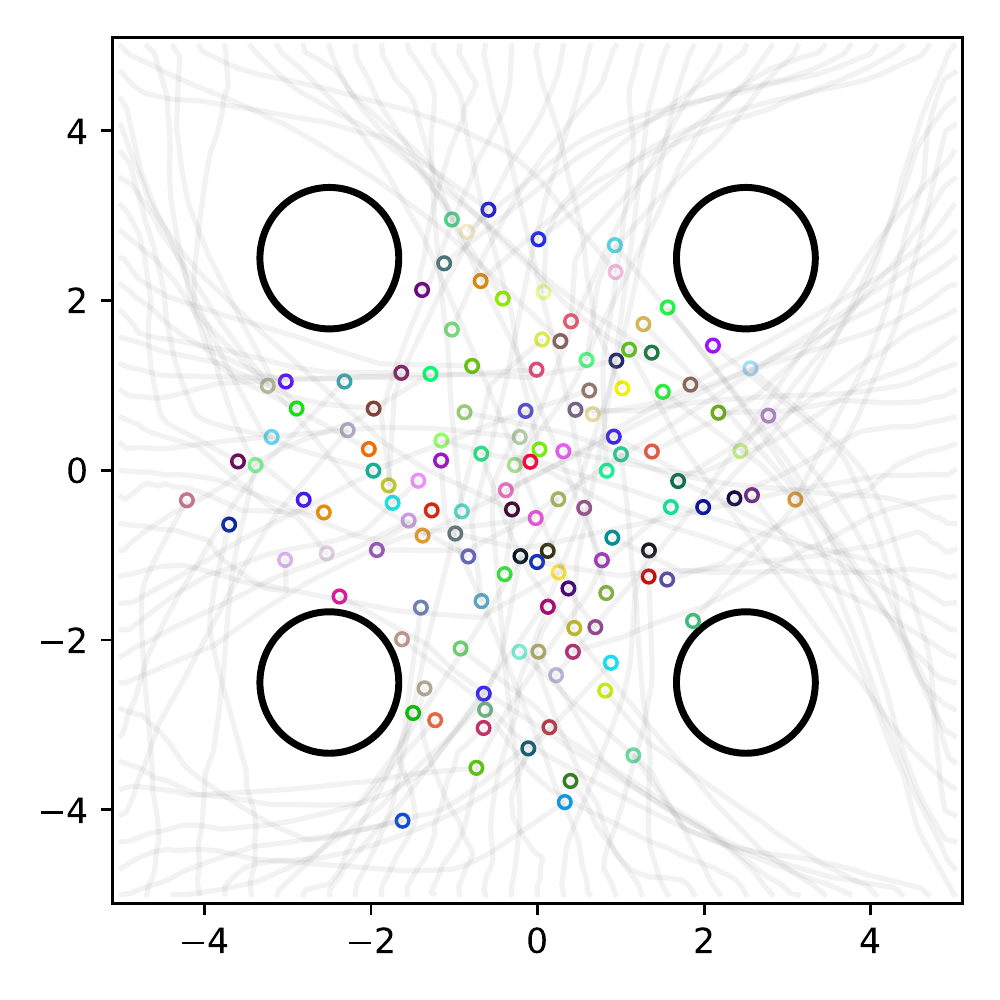}
        \caption{$t=\frac{1}{3}$}
    \end{subfigure}
    \begin{subfigure}{0.24\textwidth}
        \centering \includegraphics[width=\textwidth]{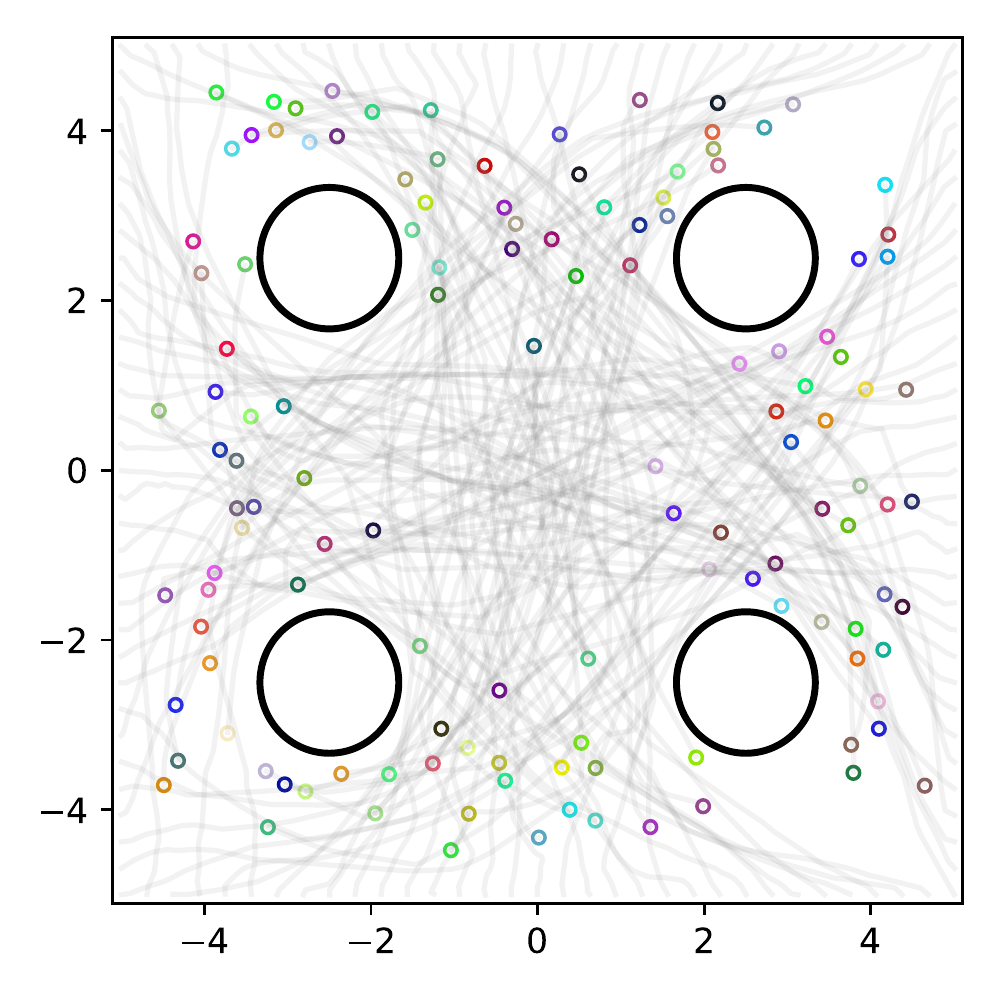}
        \caption{$t=\frac{2}{3}$}
    \end{subfigure}
    \begin{subfigure}{0.24\textwidth}
        \centering \includegraphics[width=\textwidth]{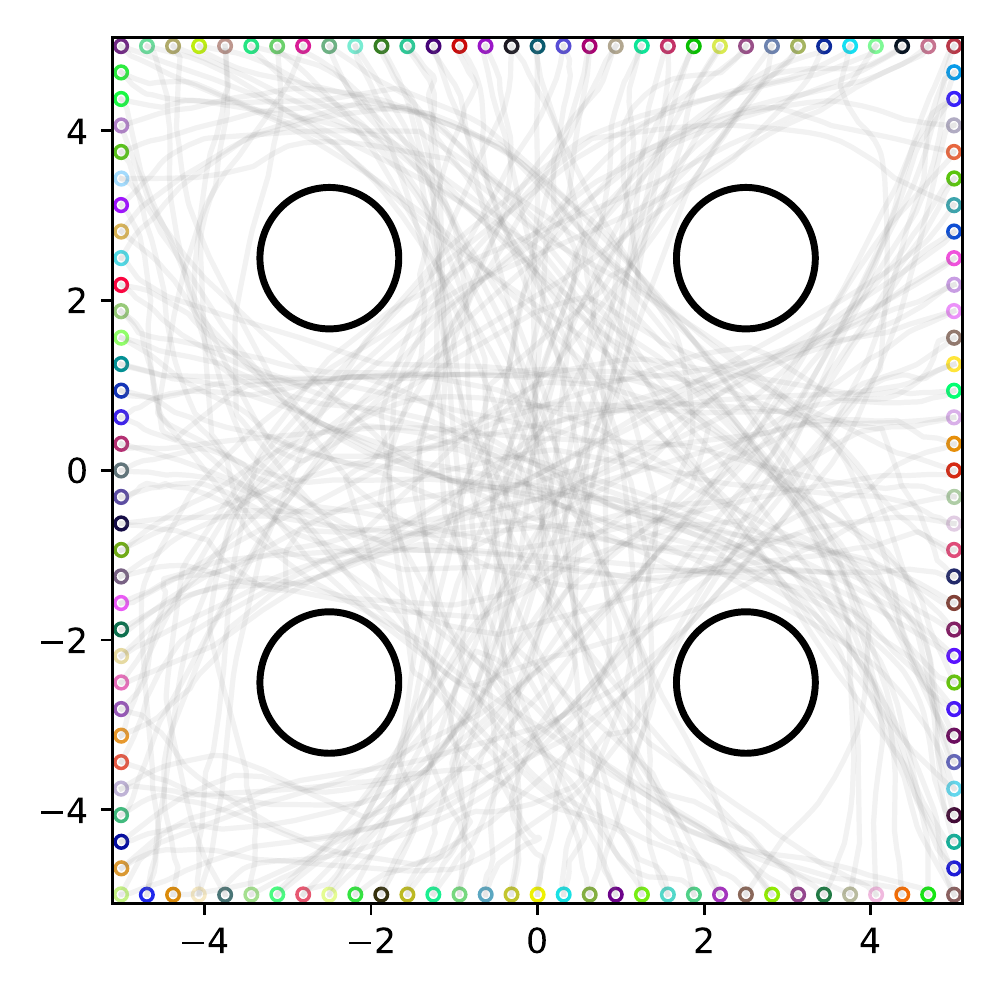}
        \caption{$t=1$}
    \end{subfigure}
    \caption{\textbf{Path planning of 128 drones with four obstacles on the plane.} We plot the predicted positions of 128 drones at time $t=0,\frac{1}{3},\frac{2}{3},1$. 
    The four black circles represent the four obstacles, and the colored circles represent the drones. 
    The paths from the initial conditions to the current positions of all drones are plotted as gray lines. There are no collisions in the three clipped snapshots, and the planned trajectories are all inside the $10\times10$ square.}
    \label{fig:test_obs_128}
\end{figure}

We also tested our SympOCnet method on this path planning problem with $128$ drones (whose radius is $0.075$) and four obstacles. The initial and terminal positions are similar to the cases of $32$ or $64$ drones (see Figure~\ref{fig:test_obs_128} (a) and (d)). The four obstacles are balls inside the room, which are represented by the black circles in Figure~\ref{fig:test_obs_128}. The constraint function $d_j$ in~\eqref{eqt:def_d} corresponding to the obstacle $E_j$ is defined by
$
    d_j(\bx) = \|\bz_j - \bx\|^2 - (C_o + C_d)^2
$
for any $\bx\in\Rn$,
where $\bz_j\in\R^2$ and $C_o>0$ are the center and radius of the obstacle $E_j$, respectively. 
The results are shown in Figure~\ref{fig:test_obs_128}. Each figure shows the current positions (represented by colored circles) and the previous trajectories (represented by grey curves) of the drones at different time $t=0,\frac{1}{3}, \frac{2}{3}, 1$. To provide a feasible solution, we set $C_d$ in the constraint function to be $0.1$ instead of $0.075$. It takes $4161$ seconds to finish $100,000$ training iterations. From the numerical results, we found that there are no collisions between obstacles and drones, and hence SympOCnet performs well and efficiently in this high-dimensional problem.

\begin{figure}[htbp]
    \centering
    \begin{subfigure}{0.24\textwidth}
        \centering \includegraphics[width=\textwidth]{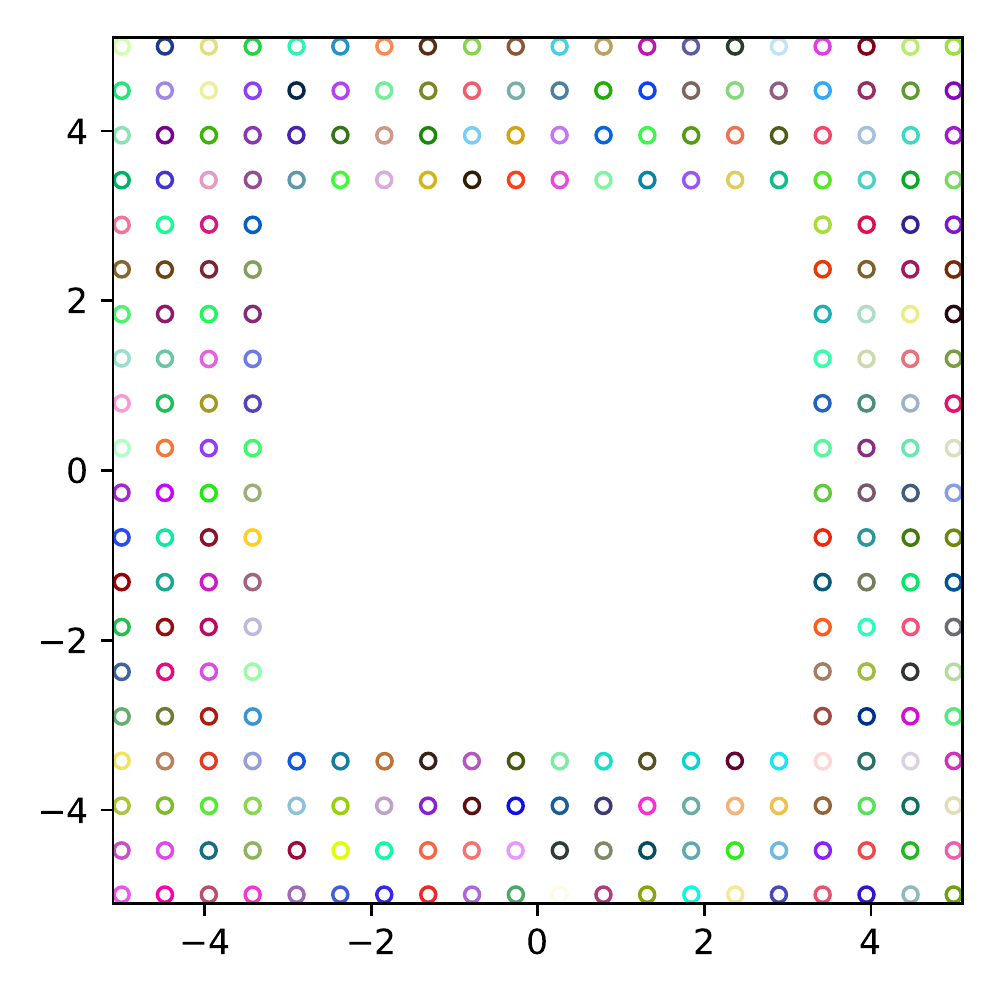}
        \caption{$t=0$}
    \end{subfigure}
    \begin{subfigure}{0.24\textwidth}
        \centering \includegraphics[width=\textwidth]{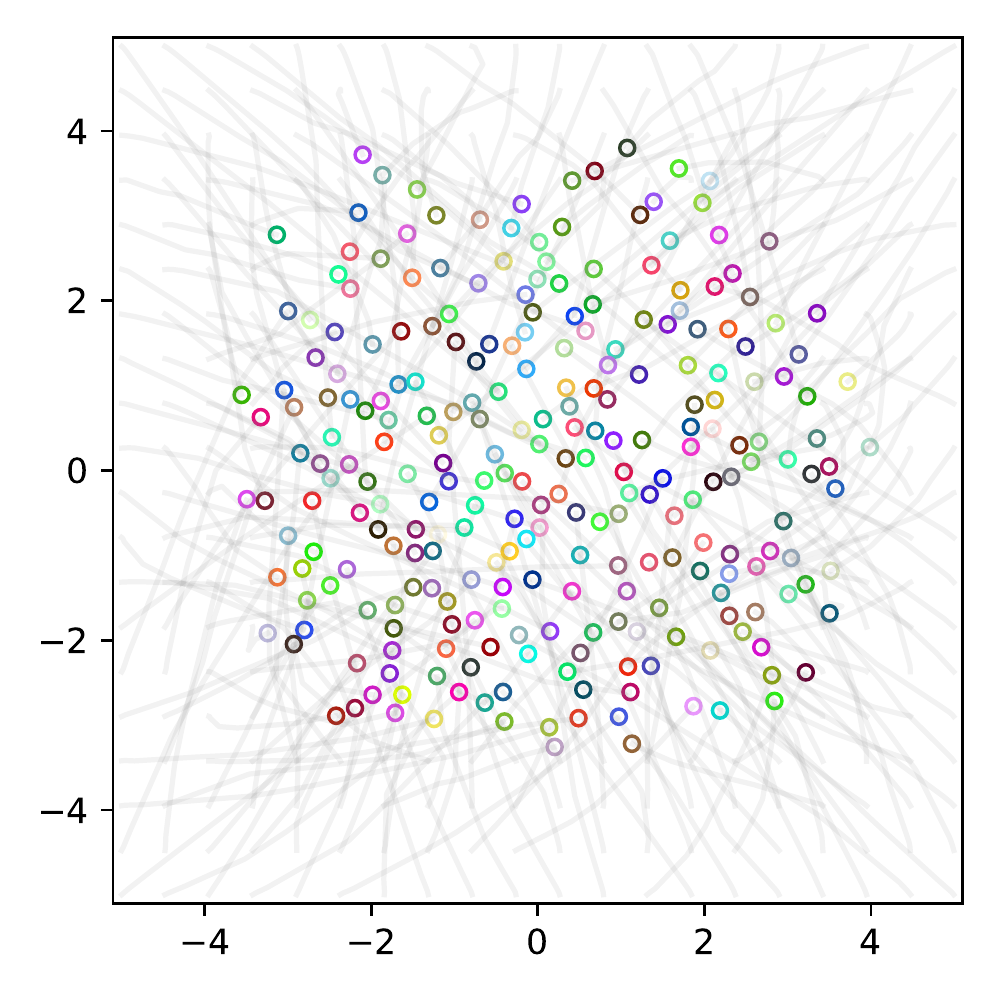}
        \caption{$t=\frac{1}{3}$}
    \end{subfigure}
    \begin{subfigure}{0.24\textwidth}
        \centering \includegraphics[width=\textwidth]{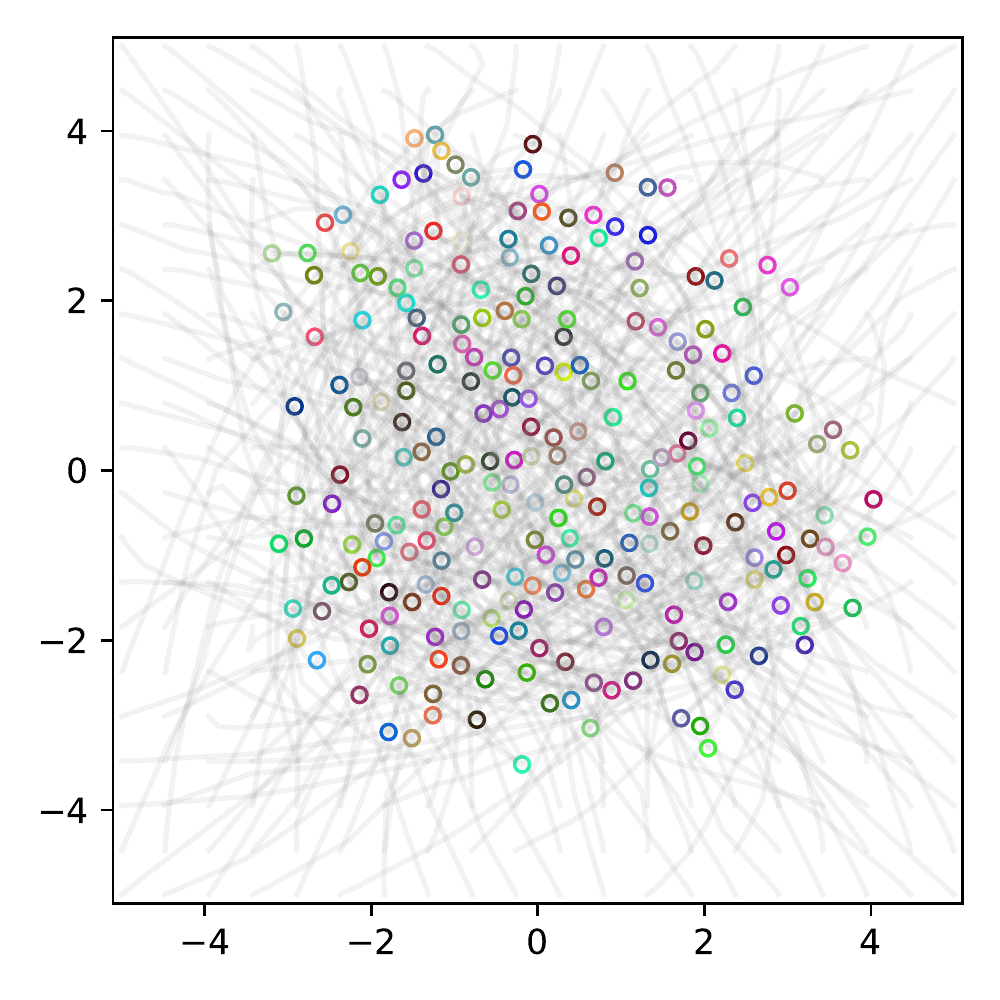}
        \caption{$t=\frac{2}{3}$}
    \end{subfigure}
    \begin{subfigure}{0.24\textwidth}
        \centering \includegraphics[width=\textwidth]{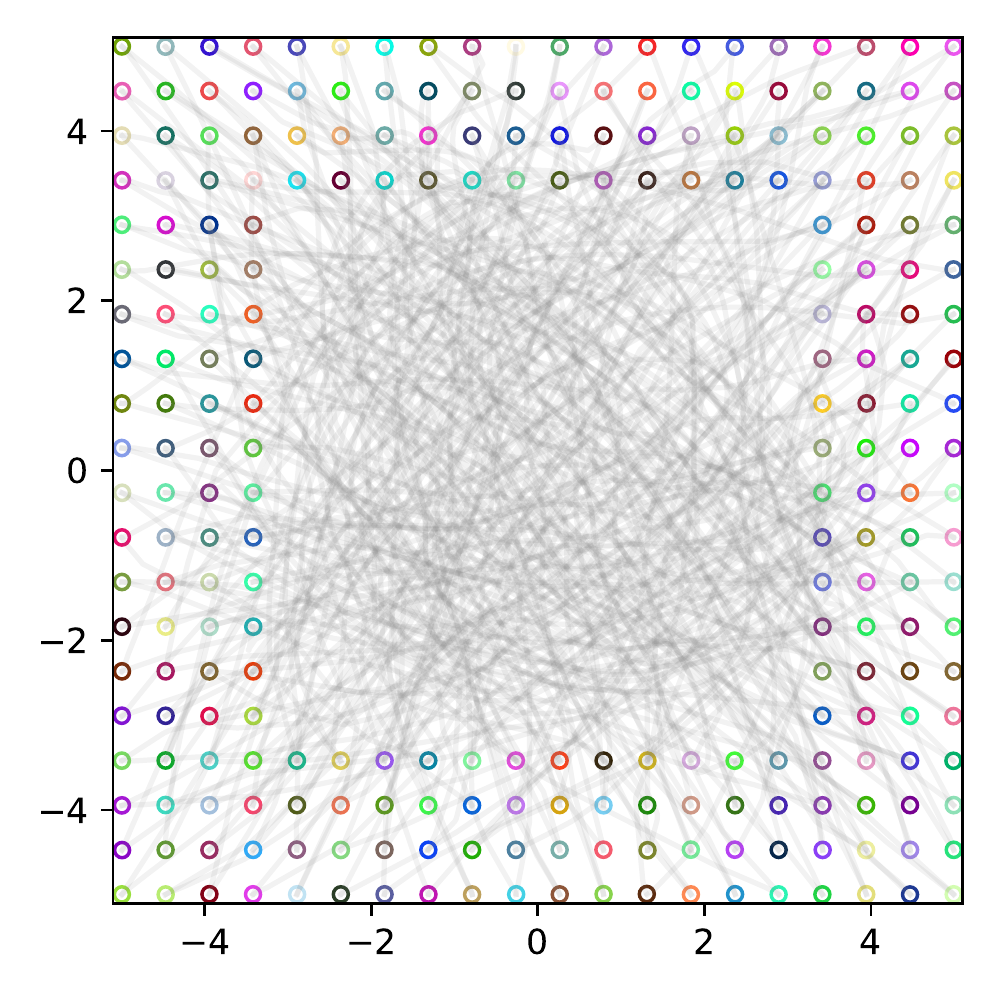}
        \caption{$t=1$}
    \end{subfigure}
    \caption{\textbf{Path planning of 256 drones on the plane.} We plot the predicted positions of 256 drones at time $t=0,\frac{1}{3},\frac{2}{3},1$. The drones are represented by colored circles. The paths from the initial positions to the current positions of all drones are plotted as gray lines.  There are no collisions in the three clipped snapshots, and the planned trajectories are inside the $10\times10$ square.}
    \label{fig:test_256}
\end{figure}

Then, we tested our proposed method on a higher dimensional problem. We set the drone number $M=256$, and hence the dimension of the state space is $n=2M=512$. We solve this path planning problem with $256$ drones (whose radius is $0.09$) in the two-dimensional room $[-5,5]^2$. The initial positions are shown in Figure~\ref{fig:test_256} (a), and the terminal condition is $\bx_T = -\bx_0$. The results are shown in Figure~\ref{fig:test_256}. To provide a feasible solution, we add a safe zone and set $C_d$ in the state constraint to be $0.1$ instead of $0.09$. From this result, we observe no collisions among the drones. It takes $5481$ seconds to obtain this feasible solution in $100,000$ training iterations. Therefore, SympOCnet is able to efficiently solve this $512$-dimensional path planning problem. 

\subsection{Multiple drones with obstacle avoidance in a three-dimensional space} \label{sec:test4_swarm}
\begin{figure}[htbp]
    \centering
    \begin{subfigure}{0.32\textwidth}
        \centering \includegraphics[width=\textwidth]{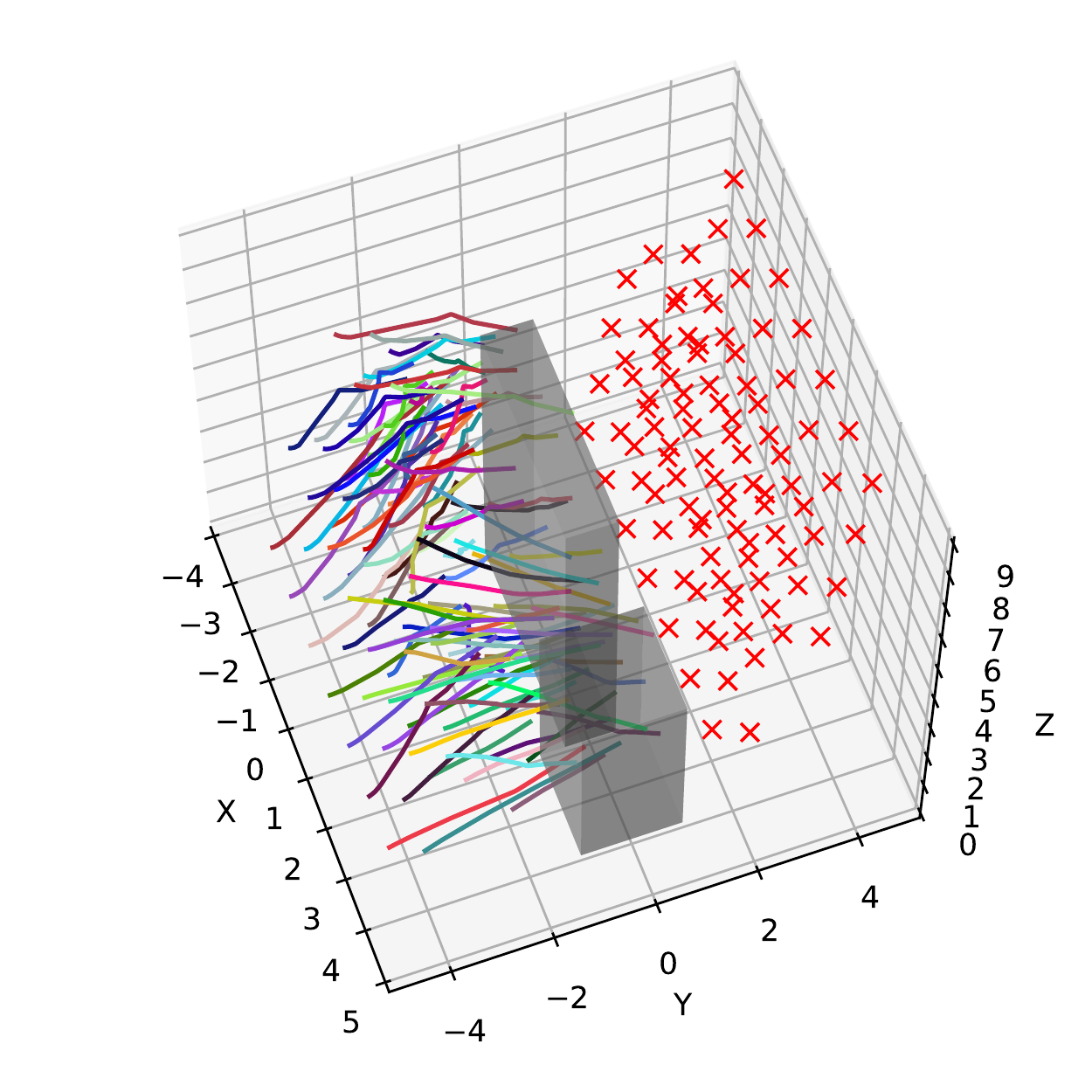}
        \caption{$t=\frac{1}{3}$}
    \end{subfigure}
    \begin{subfigure}{0.32\textwidth}
        \centering \includegraphics[width=\textwidth]{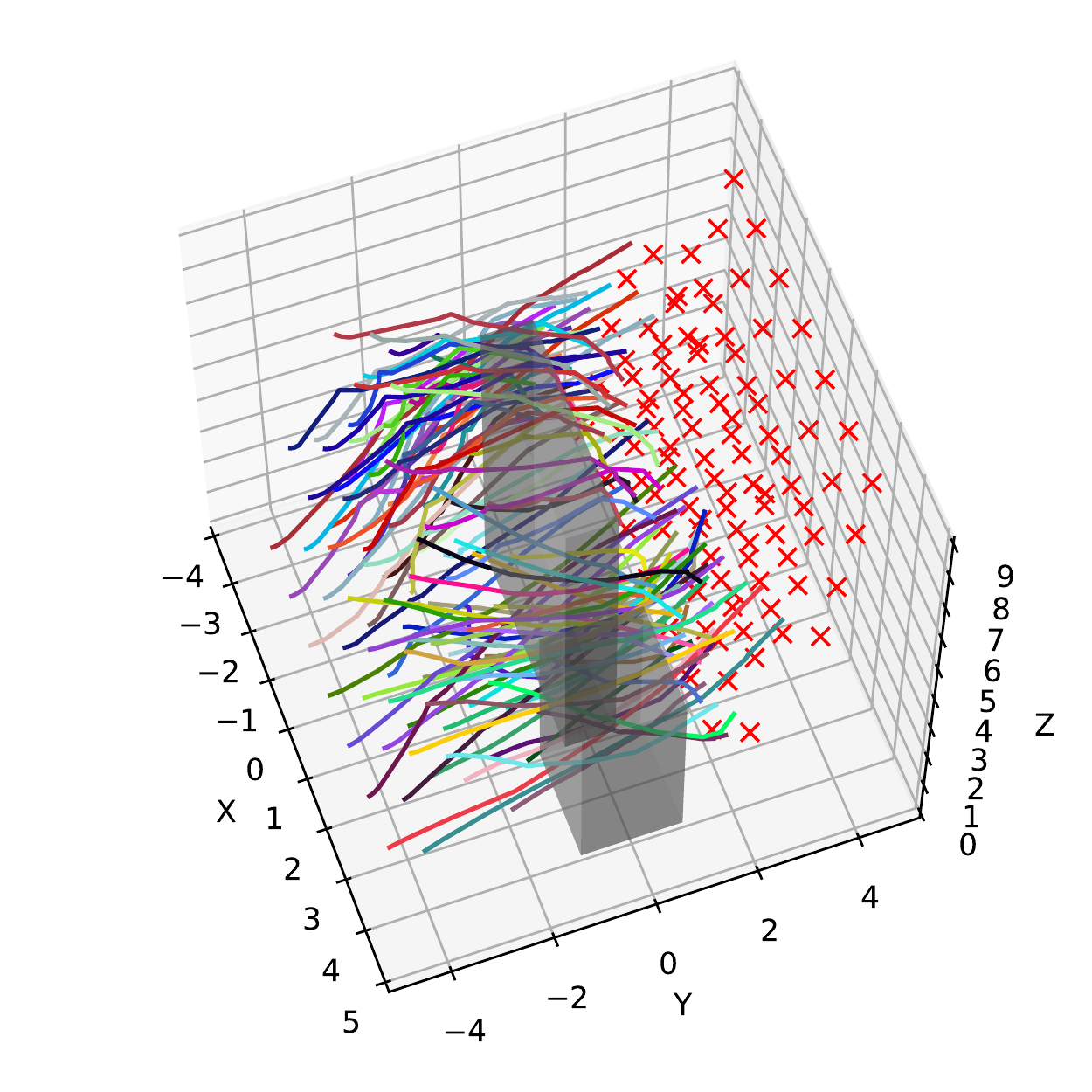}
        \caption{$t=\frac{2}{3}$}
    \end{subfigure}
    \begin{subfigure}{0.32\textwidth}
        \centering \includegraphics[width=\textwidth]{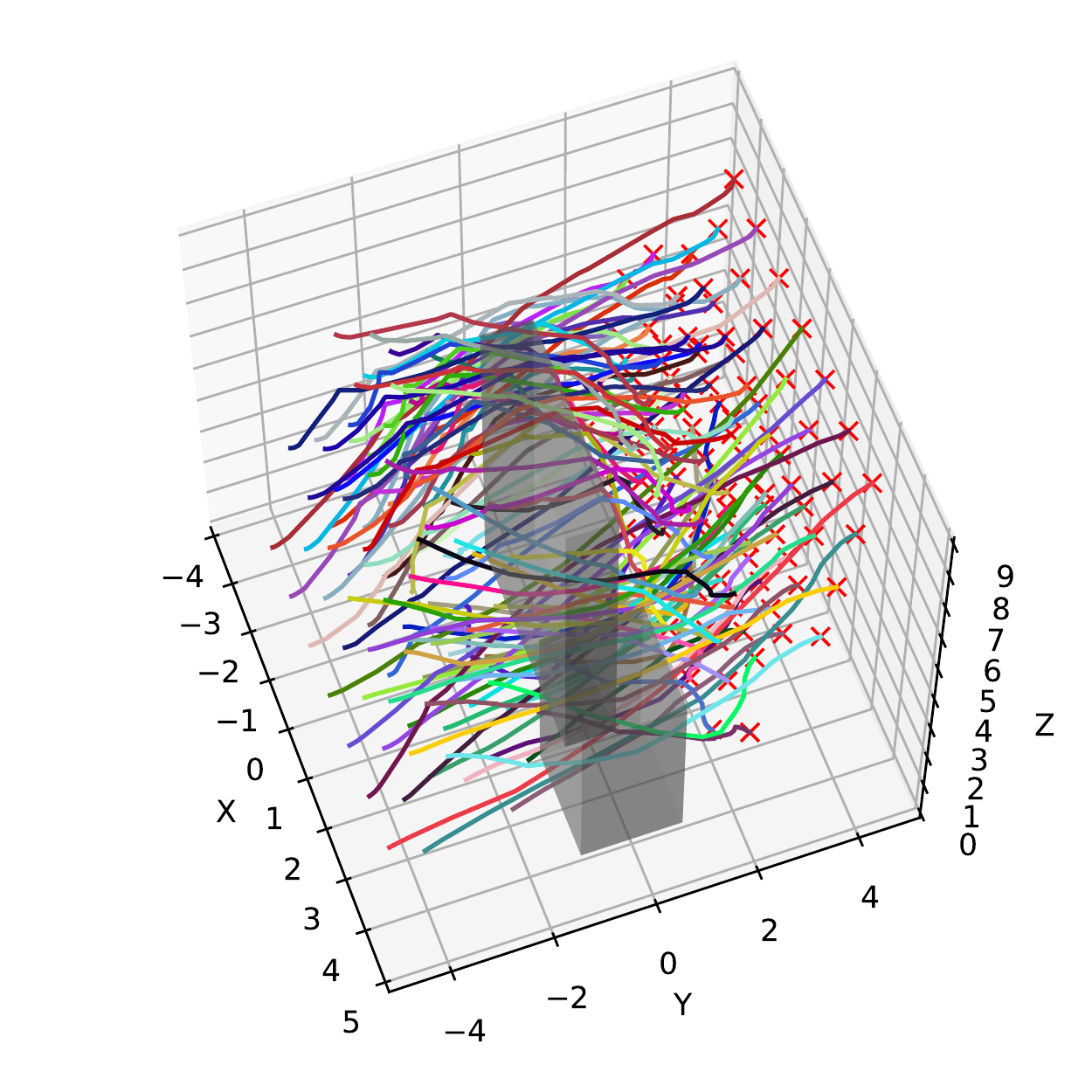}
        \caption{$t=1$}
    \end{subfigure}
    \caption{\textbf{Path planning of 100 drones in a 3d space.} We plot the predicted positions of 100 drones at time $t=\frac{1}{3},\frac{2}{3},1$. The paths from the initial positions to the current positions of all drones are plotted as colored lines. The destinations are marked as red crosses. The drones reach their destinations without any collision.}
    \label{fig:test_swarm}
\end{figure}

We consider the three-dimensional swarm path planning example in~\cite{onken2021neural}. To be specific, we consider $M = 100$ drones with radius $0.18$. Their goal is to avoid two obstacles which are located between their initial positions and terminal positions, as illustrated in Figure~\ref{fig:test_swarm}. The total dimension of the state space is $n = 3M = 300$. The $j$-th obstacle $E_j$ is defined to be the rectangular set
$[C^j_{11},C^j_{12}]\times[C^j_{21},C^j_{22}]\times[C^j_{31},C^j_{32}]$,
where $C^j_{ik}$ is a constant scalar. 
We set the constraint function $h = (h_1,h_2)$ in the same way as described before, where the function $d_j$ in~\eqref{eqt:def_d} is defined by
\begin{equation*}
\begin{split}
    &d_j(\bx) = \max\{C^j_{i1}-C_d-x_i, x_i-C^j_{i2}-C_d\colon i=1,2,3\}, 
\end{split}
\end{equation*}
for each $j=1,2$ and $\bx=(x_1,x_2,x_3)\in\R^3$. Note that $d_j(\bx)\leq 0$ if and only if $x_i\in [C^j_{i1} - C_d,C^j_{i2} + C_d]$ for each $i=1,2,3$, which holds if the ball centered at $\bx$ with radius $C_d$ intersects with the obstacle $E_j$. Therefore, the constraint $d_j(\bx)\geq 0$ enforces the drones to avoid collisions with the obstacle $E_j$.
In our experiment, we set 
\begin{equation*}
\begin{split}
    (C^1_{11}, C^1_{12},C^1_{21},C^1_{22},C^1_{31},C^1_{32}) &= (-1.8,1.8,-0.3,0.3,0.2,6.8), \\
    (C^2_{11}, C^2_{12},C^2_{21},C^2_{22},C^2_{31},C^2_{32}) &= (2.2,3.8,-0.8,0.8,0.2,3.8).
\end{split}
\end{equation*}
In other words, the constraint sets are $[-1.8,1.8]\times [-0.3,0.3]\times [0.2,6.8]$ and $[2.2,3.8]\times [-0.8,0.8]\times [0.2,3.8]$.
We apply our proposed SympOCnet method to solve this $300$-dimensional problem, and plot the result in Figure~\ref{fig:test_swarm}. Due to the high dimensionality, we do not apply the post-process. To provide a feasible solution, we set $C_d$ in the state constraint to be $0.2$ instead of $0.18$. From the numerical results, we observe no collisions among the drones and the obstacles. The minimal distance between every pair of drones is $0.3994$, which is bigger than twice of the radius. It takes $3633$ seconds to finish $150,000$ training iterations. Therefore, our proposed SympOCnet method provides a feasible solution to this high-dimensional swarm path planning problem in reasonable time.

\section{Summary}\label{sec:conclusion}
We have proposed a novel SympOCnet method for efficiently solving high-dimensional optimal control problems with state constraints.
We applied SympOCnet to several multi-agent simultaneous path planning problems with obstacle avoidance. The numerical results show SympOCnet's ability to solve high-dimensional problems efficiently with dimension more than $500$. 
These first results reveal the potential of SympOCnet for solving high-dimensional optimal control problems in real-time.
In future work, we are going to consider possible combinations with the DeepONet architecture~\cite{lu2021learning} to avoid any training cost during inference and to endow the predictive system with uncertainties, such as uncertain initial and terminal positions in path planning problems.
\section*{Acknowledgement}
The simulations were run on GeForce RTX 3090 and RTX A6000 GPU donated to us by NVIDIA.
\newpage

\bibliographystyle{siamplain}
\bibliography{biblist}

\end{document}